\documentclass[11pt,reqno]{amsart}
\usepackage{enumerate}
\usepackage{amsfonts}
\usepackage{amsmath}
\usepackage{amsthm}
\usepackage{amssymb}
\usepackage{latexsym}
\usepackage{multicol}
\usepackage{verbatim}
\usepackage{amscd}
\usepackage{hyperref}
\usepackage{amsmath, amscd}
\usepackage{pb-diagram}

\advance\textwidth by 1.2in \advance\oddsidemargin by -.6in \advance\evensidemargin by -.6in
\newtheorem*{cor}{Corollary}
\newtheorem*{lem}{Lemma}
\newtheorem*{prop}{Proposition}

\theoremstyle{definition}

\theoremstyle{definition}
\newtheorem*{thm}{Theorem}

\newtheorem*{rem}{Remark}

\newenvironment{pf}{\proof}{\endproof}
\newcounter{cnt}
\newenvironment{enumerit}{\begin{list}{{\hfill\rm(\roman{cnt})\hfill}}{%
\settowidth{\labelwidth}{{\rm(iv)}}\leftmargin=\labelwidth%
\advance\leftmargin by \labelsep\rightmargin=0pt\usecounter{cnt}}}{\end{list}} \makeatletter
\def\mydggeometry{\makeatletter\dg@YGRID=1\dg@XGRID=20\unitlength=0.003pt\makeatother}
\makeatother \theoremstyle{remark}

\numberwithin{equation}{section}

 \DeclareMathOperator{\Ht}{ht} \DeclareMathOperator{\ad}{ad}
\let\bwdg\bigwedge
\def\bigwedge{{\textstyle\bwdg}}

\begin{document}

\newcommand{\thmref}[1]{Theorem~\ref{#1}}
\newcommand{\secref}[1]{Section~\ref{#1}}
\newcommand{\lemref}[1]{Lemma~\ref{#1}}
\newcommand{\propref}[1]{Proposition~\ref{#1}}
\newcommand{\corref}[1]{Corollary~\ref{#1}}
\newcommand{\remref}[1]{Remark~\ref{#1}}
\newcommand{\defref}[1]{Definition~\ref{#1}}
\newcommand{\er}[1]{(\ref{#1})}
\newcommand{\id}{\operatorname{id}}
\newcommand{\ord}{\operatorname{\emph{ord}}}
\newcommand{\sgn}{\operatorname{sgn}}
\newcommand{\wt}{\operatorname{wt}}
\newcommand{\tensor}{\otimes}
\newcommand{\from}{\leftarrow}
\newcommand{\nc}{\newcommand}
\newcommand{\rnc}{\renewcommand}
\newcommand{\dist}{\operatorname{dist}}
\newcommand{\qbinom}[2]{\genfrac[]{0pt}0{#1}{#2}}
\nc{\cal}{\mathcal} \nc{\goth}{\mathfrak} \rnc{\bold}{\mathbf}
\renewcommand{\frak}{\mathfrak}
\newcommand{\supp}{\operatorname{supp}}
\newcommand{\Irr}{\operatorname{Irr}}
\renewcommand{\Bbb}{\mathbb}
\nc\bomega{{\mbox{\boldmath $\omega$}}} \nc\bpsi{{\mbox{\boldmath $\Psi$}}}
 \nc\balpha{{\mbox{\boldmath $\alpha$}}}
 \nc\bpi{{\mbox{\boldmath $\pi$}}}
\nc\bsigma{{\mbox{\boldmath $\sigma$}}} \nc\bcN{{\mbox{\boldmath $\cal{N}$}}} \nc\bcm{{\mbox{\boldmath $\cal{M}$}}} \nc\bLambda{{\mbox{\boldmath
$\Lambda$}}}

\newcommand{\lie}[1]{\mathfrak{#1}}
\makeatletter
\def\section{\def\@secnumfont{\mdseries}\@startsection{section}{1}%
  \z@{.7\linespacing\@plus\linespacing}{.5\linespacing}%
  {\normalfont\scshape\centering}}
\def\subsection{\def\@secnumfont{\bfseries}\@startsection{subsection}{2}%
  {\parindent}{.5\linespacing\@plus.7\linespacing}{-.5em}%
  {\normalfont\bfseries}}
\makeatother
\def\subl#1{\subsection{}\label{#1}}
 \nc{\Hom}{\operatorname{Hom}}
  \nc{\mode}{\operatorname{mod}}
\nc{\End}{\operatorname{End}} \nc{\wh}[1]{\widehat{#1}} \nc{\Ext}{\operatorname{Ext}} \nc{\ch}{\text{ch}} \nc{\ev}{\operatorname{ev}}
\nc{\Ob}{\operatorname{Ob}} \nc{\soc}{\operatorname{soc}} \nc{\rad}{\operatorname{rad}} \nc{\head}{\operatorname{head}}
\def\Im{\operatorname{Im}}
\def\gr{\operatorname{gr}}
\def\mult{\operatorname{mult}}
\def\Max{\operatorname{Max}}
\def\ann{\operatorname{Ann}}
\def\sym{\operatorname{sym}}
\def\Res{\operatorname{\br^\lambda_A}}
\def\und{\underline}
\def\Lietg{$A_k(\lie{g})(\bsigma,r)$}

 \nc{\Cal}{\cal} \nc{\Xp}[1]{X^+(#1)} \nc{\Xm}[1]{X^-(#1)}
\nc{\on}{\operatorname} \nc{\Z}{{\bold Z}} \nc{\J}{{\cal J}} \nc{\C}{{\bold C}} \nc{\Q}{{\bold Q}}
\renewcommand{\P}{{\cal P}}
\nc{\N}{{\Bbb N}} \nc\boa{\bold a} \nc\bob{\bold b} \nc\boc{\bold c} \nc\bod{\bold d} \nc\boe{\bold e} \nc\bof{\bold f} \nc\bog{\bold g}
\nc\boh{\bold h} \nc\boi{\bold i} \nc\boj{\bold j} \nc\bok{\bold k} \nc\bol{\bold l} \nc\bom{\bold m} \nc\bon{\bold n} \nc\boo{\bold o}
\nc\bop{\bold p} \nc\boq{\bold q} \nc\bor{\bold r} \nc\bos{\bold s} \nc\boT{\bold t} \nc\boF{\bold F} \nc\bou{\bold u} \nc\bov{\bold v}
\nc\bow{\bold w} \nc\boz{\bold z} \nc\boy{\bold y} \nc\ba{\bold A} \nc\bb{\bold B} \nc\bc{\bold C} \nc\bd{\bold D} \nc\be{\bold E} \nc\bg{\bold
G} \nc\bh{\bold H} \nc\bi{\bold I} \nc\bj{\bold J} \nc\bk{\bold K} \nc\bl{\bold L} \nc\bm{\bold M} \nc\bn{\bold N} \nc\bo{\bold O} \nc\bp{\bold
P} \nc\bq{\bold Q} \nc\br{\bold R} \nc\bs{\bold S} \nc\bt{\bold T} \nc\bu{\bold U} \nc\bv{\bold V} \nc\bw{\bold W} \nc\bz{\bold Z} \nc\bx{\bold
x} \nc\KR{\bold{KR}} \nc\rk{\bold{rk}} \nc\het{\text{ht }}

\nc\toa{\tilde a} \nc\tob{\tilde b} \nc\toc{\tilde c} \nc\tod{\tilde d} \nc\toe{\tilde e} \nc\tof{\tilde f} \nc\tog{\tilde g} \nc\toh{\tilde h}
\nc\toi{\tilde i} \nc\toj{\tilde j} \nc\tok{\tilde k} \nc\tol{\tilde l} \nc\tom{\tilde m} \nc\ton{\tilde n} \nc\too{\tilde o} \nc\toq{\tilde q}
\nc\tor{\tilde r} \nc\tos{\tilde s} \nc\toT{\tilde t} \nc\tou{\tilde u} \nc\tov{\tilde v} \nc\tow{\tilde w} \nc\toz{\tilde z}

\title{A categorical approach to Weyl modules}

\author{Vyjayanthi Chari}
\address{Department of Mathematics, University of California, Riverside, CA 92521, USA}
\email{chari@math.ucr.edu}
\author{Ghislain Fourier}
\address{Mathematisches Institut, Universit\"at zu K\"oln, Germany}
\email{gfourier@mi.uni-koeln.de}
\author{Tanusree Khandai}
\address{Harish-Chandra Research Institute, Allahabad, India}
\email{tanusree@mri.ernet.in}
\thanks{VC was partially supported by the NSF grant DMS-0500751}
\thanks{G.F. was supported by the DFG-project \lq\lq Kombinatorische Beschreibung von Macdonald und Kostka-Foulkes Polynomen \rq\rq}
\begin{abstract} Global and local  Weyl Modules were introduced via generators and relations in the context of affine Lie algebras in \cite{CPweyl} and were motivated by representations of quantum affine algebras. In \cite{FL} a more general case was considered by replacing the polynomial ring with the coordinate ring of an algebraic variety and partial results analogous to those in \cite{CPweyl} were obtained.
  In this paper, we show that there is a natural definition of the local and global Weyl modules via homological properties. This characterization allows us to define the Weyl functor from the category of left modules of a commutative algebra to the category of modules for a simple Lie algebra. As an application we are able to understand the relationships of these functors to tensor products, generalizing results in \cite{CPweyl} and \cite{FL}.  We also analyze the fundamental Weyl modules and show that unlike the case of the affine Lie algebras,  the Weyl   functors  need not be left exact.
\end{abstract}
\maketitle

\section{Introduction}
The category of finite--dimensional representations of affine and quantum affine Lie algebras has been intensively studied in recent years. One of the reasons that this category has proved to be interesting is the fact that it is not semi-simple. Moreover, it was proved in \cite{CPweyl}  that irreducible representations of the quantum affine algebra specialized to reducible indecomposable  representations of the affine Lie algebra. This phenomenon is analogous to the one observed in modular representation theory where an irreducible finite--dimensional  representation in characteristic zero becomes reducible on passing to characteristic $p$ and is called a Weyl module.

The definition of Weyl modules (global and local) in \cite{CPweyl} for affine algebras was motivated by this analogy. Thus given any dominant integral weight of the semisimple Lie algebra $\lie g$, one can define an infinite--dimensional left module  $W(\lambda)$ for the corresponding affine (in fact for the loop) algebra via generators and relations. The module $W(\lambda)$ is a direct sum of finite--dimensional $\lie g$--modules and it was shown in \cite{CPweyl} that it is also a right module for a polynomial algebra $\mathbb A_\lambda$ which is canonically associated with $\lambda$. The local Weyl modules are obtained by tensoring the global Weyl modules with irreducible modules for $\mathbb A_\lambda$ or equivalently can be given via generators and relations. A necessary and sufficient condition for the tensor product of local Weyl modules to be a local Weyl module was given. Using this fact, the  character of the local Weyl module  was conjectured in \cite{CPweyl} and the conjecture was heavily influenced by the connection with quantum affine algebras. In particular, the conjecture implied that the dimension of the local Weyl module was independent of  the choice of the irreducible $\mathbb A_\lambda$--module, i.e that the global Weyl module is a free module for $\mathbb A_\lambda$.
The character formula was proved in \cite{CPweyl} for $\lie {sl_2}$, in   \cite{CL} for $\lie{sl}_{r+1}$, in \cite{FoL} for simply--laced algebras and the general case can be deduced by passing to the quantum case by using the work of \cite{K} and \cite{BN}.

In \cite{FL}, Feigin and Loktev extended the notion of Weyl modules to the higher--dimensional case, i.e. instead of the  loop algebra they worked with the Lie  algebra $\lie g\otimes A$ where $A$ is the coordinate ring of an algebraic variety and obtained analogs of some of the results of \cite{CPweyl}. For instance when $\lie g$ is of type $\lie{sl}_2$ and  $A$ is the polynomial ring in two variables they compute the dimension of the Weyl module. They also  give a necessary and sufficient condition for the tensor product of local Weyl modules to be a local Weyl module analogous to the one in \cite{CPweyl}.  However, they do not define the algebra $\mathbb A_\lambda$ and the bi--module structure on $W(\lambda)$ and hence do not say much about the structure of the global Weyl module.

In this paper, we take a more general functorial approach to Weyl modules associated to the algebra $\lie g\otimes A$, where $A$ is a commutative associative algebra (with unit) over the complex numbers. This approach  (as also the approach in \cite{CG1}, \cite{CG2}) is motivated by the methods used to study another well--known category in representation theory: the  BGG-category $\cal O$ for semi--simple Lie algebras. As a result we are able to extend the definition of Weyl modules to a more general situation and allows us to do a deeper analysis of the global Weyl modules. We also give the classification and description of irreducible modules for $\lie g\otimes A$ for an arbitrary finitely generated algebra which is analogous to the one given in  \cite{C1},\cite{CP1},\cite{L},\cite{R} in the case when $A$ is a polynomial algebra.

We now explain our results in some detail. Let $\cal I_A$ be the category of   $\lie g\otimes A$--modules which are integrable as $\lie g$--modules.  For $\lambda\in P^+$ we let $\cal I^\lambda_A$ be the full subcategory of $\cal I_A$ consisting of objects whose weights are bounded above by $\lambda$.
Given $\lambda\in P^+$, one can define in a canonical way a  projective module $P_A(\lambda)\in\cal I_A$ and  we prove that the global Weyl module $W_A(\lambda)$ is the largest quotient  of $P_A(\lambda)$ that lies in $\cal I^\lambda_A$.  We then define a right action of the algebra $\bu(\lie h\otimes A)$ on $W_A(\lambda)$ where $\lie h$ is a Cartan subalgebra of $\lie g$ which is compatible with the left action of $\lie g\otimes A$. Let $\ba_\lambda$ be the quotient of $\bu(\lie h\otimes A)$ by the torsion ideal for this action so that $W_A(\lambda)$ can be regarded as a bi-module for $(\lie g\otimes A, \ba_\lambda)$. We prove that the bimodule structure is functorial in $A$.

Let $\bw^\lambda_A$ be the right exact functor $W_A(\lambda) \otimes_{\ba_\lambda}$ from the category $\mode\ba_\lambda$ of left modules for $\ba_\lambda$ to $\cal I^\lambda_A$. The local Weyl modules are then just $\bw^\lambda_AM$ where $M$ is an irreducible object of $\mode\ba_\lambda$. In section 3, we prove that one can define a functor $\br^\lambda_A$ which is exact and right adjoint to $\bw^\lambda_A$. That allows us to give a categorical characterization of the local Weyl modules and more generally of  the modules $\bw_A^\lambda M$, $M\in\mode\ba_\lambda$.  Namely we prove that these modules are given by the vanishing of $\Hom_{\cal I^\lambda_A}$ and $\Ext^1_{\cal I^\lambda_A}$ and we show also that the functors $\bw^\lambda_A$ are left exact iff we have vanishing of $\Ext^2_{\cal I^\lambda_A}$.

In section 4 we prove that the algebra $\ba_\lambda$ is finitely generated iff $A$ is finitely generated. We  use the results of section 3 to study the relationship between the functors $\bw^{\lambda+\mu}_{A\oplus B}$ and $\bw^\lambda_A\otimes \bw^\mu_B$ when $A,B$ are finite--dimensional algebras. In section 5, we give a necessary and sufficient condition for the tensor product $\bw_A^\lambda M\otimes\bw_A^\mu N$ to be isomorphic to $\bw^{\lambda+\mu}_A(M\otimes N)$ when $A$ is finitely generated and $M,N\in\mode\ba_\lambda$.

In section 6 we assume that $A$ is finitely generated and that the Jacobson radical of $A$ is $0$. We prove that the algebra $\ba_\lambda$ is isomorphic to the ring of invariants of a subgroup $S_\lambda$ of the symmetric group on $d_\lambda$ letters acting on $A^{\otimes d_\lambda}$.  Here $d_\lambda$ is a positive integer naturally associated with $\lambda$. This implies that  the irreducible modules in $\mode\ba_\lambda$ are determined (up to isomorphism) by the orbits of this action.

The tensor product results  of Sections 4 and 5 imply that to understand the local Weyl modules  it is enough to understand local Weyl modules corresponding to certain special orbits. In section 7, we consider the case when $\xi$ is the orbit of a point in $A^{\otimes d_\lambda}$   which has trivial stabilizer under the entire symmetric group  $S_{d_\lambda}$. In this case $\bw^\lambda_AM_\xi$ is a tensor product of the  local fundamental Weyl modules and we describe the character of these modules completely for any finitely generated algebra $A$ and for the classical simple Lie algebras.

The results of section 7 show that there are many important differences between the study of Weyl modules for the polynomial algebra in one variable and the more general case considered here.  The dimension of the local fundamental Weyl modules associated to $A$ depends on $\xi$ if the variety associated to $A$ is not smooth. It also proves that the dimension of $\bw_A^\lambda M_\xi$ is not independent of $\xi$ even if $A$ is an irreducible smooth variety and $\xi$ is the orbit of a point in $A^{\otimes d_\lambda}$  with trivial stabilizer for the $S_{\lambda}$-action. In particular, this proves that the global Weyl module is not projective as a right $\ba_{\lambda}$--module (and hence the Weyl functors not exact) even when $A$ is the polynomial ring in two variables. There are thus,  many natural and interesting algebraic and geometric  questions that  arise as a result of this paper which will be studied elsewhere.

\vskip12pt
{\em Acknowledgements: We would like to thank Wee Liang Gan, Michael Ehrig, Friederich Knop, Peter Littelmann for many discussions on the algebra $\mathbb A_\lambda$.  We are grateful to Peter Russell for his patience with our long discussions and our not always well-formulated questions  on group actions, homological algebra and commutative algebra. Finally, particular thanks are due to Shrawan Kumar for sharing with us, his result (Proposition \ref{kunneth} ) on extensions between tensor products of modules for direct sums of Lie algebras.}

\section{Preliminaries }

\subsection{}
 Throughout the paper $\bc$  denotes the set of complex  numbers and
$\bz_+$ the set of non--negative integers.  Let  $\lie g$   be a finite--dimensional
   simple  Lie algebra of rank $n$ with Cartan matrix $(a_{ij})_{i,j\in I}$ where $I=\{1,\cdots, n\}$. Fix a Cartan subalgebra $\lie h$ of $\lie g$  and let
  $R$ denote the corresponding  set of   roots.  Let $\{\alpha_i\}_{i\in I}$ (resp. $\{\omega_i\}_{i\in I}$)  be  a
set of simple roots (resp. fundamental weights) and  $Q$ (resp. $Q^+$), $P$ (resp. $P^+$) be the integer span (resp. $\bz_+$--span) of the
simple roots and fundamental weights respectively. Denote by  $\le $  the usual partial order on $P$,  $$\lambda,\mu\in
P,\ \ \lambda\le \mu\ \iff\  \mu-\lambda\in Q^+.$$ Set $R^+= R\cap Q^+$ and let  $\theta$ be  the unique maximal element in $R^+$ with
respect to the partial order.

 Let $x^\pm_\alpha$, $h_i$, $\alpha\in R^+$, $i\in I$ be a Chevalley basis of $\lie g$ and set $x_i^\pm=x^\pm_{\alpha_i}$,
 $h_\alpha=[x^+_\alpha, x^-_\alpha]$ and note that $h_i=h_{\alpha_i}$. For each $\alpha\in R^+$, the  subalgebra of $\lie g$ spanned by $\{x^\pm_\alpha, h_\alpha\}$ is isomorphic to $\lie{sl}_2$.
 Define subalgebras  $\lie n ^\pm $ of $\lie g$,
  by  $$\lie n^\pm=\bigoplus_{\alpha\in R^+}\bc x^\pm_\alpha,$$
   and note that $$\lie g=\lie n^-\oplus\lie h\oplus\lie n^+.$$ Given any Lie algebra $\lie a$,
    let $\bu(\lie a)$ be the universal enveloping algebra of $\lie a$.  The map
    $x\to x\otimes 1+1 \otimes x$, $x\in\lie a$ extends to an algebra homomorphism $\Delta: \bu(\lie a)\to \bu(\lie a)\otimes\bu(\lie a)$.
     By the Poincare Birkhoff Witt theorem, we know that if $\lie b$ and $\lie c$ are Lie subalgebras of $\lie a$ such that   $\lie a=\lie b\oplus\lie c$ as vector spaces then $$\bu(\lie a)\cong\bu(\lie b)\otimes\bu(\lie c)$$ as vector spaces.

\subsection{} 
Let $A$ be a commutative associative algebra with unity over $\bc$  and let $A_+$ be  a fixed vector space complement to the subspace $\bc$ of $A$.  Given a Lie algebra $\lie a$  define a Lie algebra structure on $\lie a\otimes A$,
 by $$[x\otimes a, y\otimes b]=[x, y]\otimes ab,\ \ x,y\in\lie g, \ \ a,b\in A.$$
 If  $\phi:B\to A$ is a homomorphism of associative  algebras, there exists a corresponding homomorphism $\phi_{\lie a}:\lie a\otimes B\to\lie a\otimes A$
 of Lie algebras, which is injective (resp. surjective) if $\phi$ is injective (resp. surjective).
 In particular,  if $B$ is a subalgebra of $A$, the Lie algebra $\lie a\otimes B$
  can be regarded naturally as a Lie subalgebra of $\lie a\otimes A$
 and we identify  $\lie a$ with the  Lie subalgebra  $\lie a\otimes \bc$ of $\lie a \otimes A$.
   Similarly, if $\lie b$ is a Lie subalgebra of $\lie a$, then $\lie b\otimes A$ is
    naturally isomorphic to a subalgebra of $\lie a\otimes A$. Finally we denote by $\bu(\lie g\otimes A_+)$ the subspace of $\bu(\lie g\otimes A)$ spanned by monomials in the elements $x\otimes a$ where $x\in \lie g$, $a\in A_+$. The following is elementary but we include a proof for the reader convenience and because it is used repeatedly throughout the paper.
    \begin{lem}\label{ideal} Let $\lie g$ be a finite--dimensional simple Lie algebra and $A$ a commutative associative algebra with unity over $\bc$. Then any ideal of $\lie g\otimes A$ is of the form $\lie g\otimes S$ for some ideal $S$ of $A$ and $[\lie g\otimes A/S,\lie g\otimes A/S]=\lie g\otimes A/S.$
    \end{lem}
\begin{pf} Let $\lie i$ be  an ideal in $\lie g\otimes A$  and set $$S=\{a\in A:\lie g\otimes a\subset\lie i\}.$$ Since $\lie g=[\lie g,\lie g]$ we see that $S$ is an ideal on $A$.  The Lemma follows if we prove that $\lie g\otimes S=\lie i.$ Let $x\in\lie i$ and write
 $$x=\sum_{\alpha\in R}x_\alpha\otimes  a_\alpha+\sum_{i\in I} h_i\otimes a_i,$$ for some $a_\alpha, a_i\in A$. We proceed by induction on $$r=\#\{\alpha\in R:a_\alpha\ne 0\},$$ to show that $\lie g\otimes a_\alpha\subset \lie i$ and $\lie g\otimes a_i\subset\lie i$ for all $\alpha\in R$, $i\in I$. If $r=0$, we have $$[\sum_{i\in I} h_i\otimes a_i,x_j^+]=x_j^+\otimes \sum_{i\in I}\alpha_j(h_i)a_i\in\lie i,\ \ j\in I.$$ Since the Cartan matrix of $A$ is invertible, it follows now that $x_j^+\otimes a_i\in\lie i$ for all $i,j\in I$ and since $\lie g$ is simple we see that $\lie g\otimes a_i\in\lie i$ for all $i\in I$.

Suppose now that we have proved the result when $0\le r<k$ and suppose that $a_{\beta_1},\cdots,a_{\beta_k}$ are the non--zero elements. Choose $h\in\lie h$ such that $\beta_k(h)\ne 0$ and $\beta_{k-1}(h)=0$. Then
$$0\ne [h,x]=\sum_{s=1}^{k-2}\beta_s(h)x_\alpha\otimes a_{\beta_s}+ \beta_k(h)x_{\beta_k}\otimes a_{\beta_k}\in\lie i.$$ The induction hypothesis applies to $[h,x]$ and we find that $$a_{\beta_{k}}\in S,\ \ \ x-x_{\beta_{k}}\otimes a_{\beta_{k}}\in\lie i.$$ The induction hypothesis again applies to $x-(x_{\beta_k}\otimes a_{\beta_k})$  and we get the result.
\end{pf}

\subsection{}
Let $V$ be any $\lie g$--module. We say that $V$ is locally finite--dimensional if
any element  of $ V$ lies in a finite--dimensional $\lie g$--submodule of $V$.  This means that $V$ is isomorphic to a direct
sum of irreducible finite--dimensional $\lie g$--modules and hence we can write $$V=\bigoplus_{\lambda\in\lie h^*} V_\lambda,$$ where
$V_\lambda=\{v\in V:hv=\lambda(h)v,\ \ \forall\ h\in\lie h\}$. We set
$$\wt(V)=\{\lambda\in\lie h^*:V_\lambda\ne 0\}.$$
 For $\lambda\in P^+$, let $V(\lambda)$ be the simple
$\lie g$--module which is  generated by an element $v_\lambda\in V(\lambda)$ satisfying the defining relations:
$$
\lie n^+ v_\lambda=0,\quad hv_\lambda=\lambda(h)v_\lambda,\quad (x^-_{i})^{\lambda(h_i)+1}v_\lambda =0,
$$
for all~$h\in\lie h$, $i\in I$. Then, $$\wt(V(\lambda))\subset \lambda -Q^+,\ \ \dim V(\lambda)<\infty.$$ Moreover  any irreducible locally
finite--dimensional $\lie g$--module is isomorphic to $V(\lambda)$ for some $\lambda\in P^+$. The following can be found in \cite{B}.
\begin{lem}\label{sscond} Let $\lie a$ be a Lie algebra such that $[\lie a,\lie a]=\lie a$ and assume that $\lie a$ has a faithful  finite--dimensional irreducible representation. Then $\lie a$ is a semi--simple Lie algebra.
\end{lem}

\subsection{}
Suppose that $\lie g$ is a finite--dimensional semisimple Lie algebra and that $\lie g_1$, $\lie g_2$ are ideals of $\lie g$ such that $$\lie g\cong \lie g_1\oplus\lie g_2$$ as Lie algebras. Then $\lie g_1$ and $\lie g_2$ are also semisimple Lie algebras and it is standard that any irreducible  finite--dimensional representation of ${\lie g}$ is isomorphic to a tensor product of irreducible representations of $\lie g_1$ and $\lie g_2$.

 \begin{prop}\label{ideal-irred} Let $A$ and $B$ be  commutative associative algebras. Any finite--dimensional  irreducible representation $V$ of $\lie g\otimes (A\oplus B)$ is isomorphic to a tensor product $V_1\otimes V_2$  where $V_1$ and $V_2$ are irreducible representations of $\lie g\otimes A$ and $\lie g\otimes B$ respectively.
\end{prop}
\begin{pf} Let $\rho:\lie g\otimes(A\oplus B)\to\End(V)$ be
 an irreducible finite--dimensional representation. Then
 $\ker\rho$ is an ideal of finite codimension in $\lie g\otimes (A\oplus B)$
 and hence $$\ker\rho= \lie g\otimes M,$$ for some ideal $M$ of $A\oplus B$.
  Since any ideal of $A\oplus B$ is of the form $M_1\oplus M_2$ where $M_1,M_2$ are ideals in $A$ and $B$ respectively,
   we see that $V$ is a faithful irreducible  representation of $\tilde{\lie g}=\lie g\otimes (A/M_1\oplus B/M_2)$.
 Lemma \ref{sscond} implies  that $\tilde{\lie g}$ is a finite--dimensional semi-simple Lie algebra. The result
 now follows by the comments preceding the statement of this  proposition.
\end{pf}

\subsection{}
We shall need  the following result due to Shrawan Kumar \cite{Ku}.
\begin{prop}\label{kunneth}
 For $r=1,2$, let $\lie g_r$  be a finite--dimensional Lie algebra and assume that $U_r, V_r$ are  finite dimensional $\lie g_r$--modules. For all $m \geq 0$ we have
$$ \Ext^{m}_{\lie g_1 \oplus \lie g_2}(U_1 \otimes U_2, V_1 \otimes V_2) \cong \bigoplus_{p +q = m} \Ext^{p}_{\lie g_1}(U_1, V_1) \otimes \Ext^{q}_{\lie g_2}(U_2, V_2) .$$
\end{prop}

\section{The category $\mathcal{I}_A $}

\subsection{}
 Let $\cal I_A$   be the category whose
objects are modules for $\lie g\otimes A$ which are locally finite--dimensional $\lie g$--modules and morphisms $$\Hom_{\cal
I_A}(V,V')=\Hom_{\lie g\otimes A}(V,V'),\ \ V,V'\in\cal I_A.$$ Clearly  $\cal I_A$ is an abelian category and is closed under tensor products.
 We shall use the following elementary result often without  mention in the rest of the paper.
 \begin{lem}\label{superelem} Let $V\in\Ob\cal I_A$.  \begin{enumerit} \item If $V_\lambda\ne 0$ and  $\wt V\subset\lambda-Q^+$, then $\lambda\in P^+$ and $$(\lie n^+\otimes A) V_\lambda=0,\ \ (x_i^-)^{\lambda(h_i)+1}V_\lambda=0,\ \ i\in I.$$  If in addition, $V=\bu(\lie g\otimes A)V_\lambda$ and $\dim V_\lambda=1$, then $V$ has a unique irreducible  quotient.
 \item If $V=\bu(\lie g\otimes A)V_\lambda$ and $(\lie n^+\otimes A)V_\lambda=0$, then $\wt(V)\subset\lambda-Q^+$.
 \item If $V\in\cal I_A$ is irreducible and finite--dimensional, then there exists $\lambda\in\wt V$ such that $$ \dim V_\lambda=1,\ \ \wt(V)\subset\lambda- Q^+.
$$
     \hfill\qedsymbol
 \end{enumerit}
 \end{lem}

\subsection{}
Regard $\bu(\lie g\otimes A)$ as a right $\lie g$--module via right multiplication and
given a  left $\lie g$--module $V$, set $$P_A(V)=\bu(\lie g\otimes A)\otimes_{\bu(\lie g)} V. $$
  Then  $P_A(V)$ is a left $\lie g\otimes A$--module by left multiplication
  and we have an isomorphism of vector spaces \begin{equation}\label{isov} P_A(V)\cong \bu(\lie g\otimes A_+)\otimes_\bc V.\end{equation}

  \begin{prop}\label{projelem} Let $V$ be a locally finite--dimensional $\lie g$--module.
  Then $P_A(V)$ is a projective object of $\cal I_A$.  If in addition $V\in\cal I_A$, then the map $P_A(V)\to V$ given by $u\otimes v\to uv$
  is a surjective morphism of objects in $\cal I_A$. Finally, if $\lambda\in P^+$, then $P_A(V(\lambda))$  is generated as a $\bu(\lie g\otimes A)$--module by the element $p_\lambda=1\otimes v_\lambda$ with defining relations \begin{equation}\label{defproj}
\lie n^+ p_\lambda=0,\quad h p_\lambda=\lambda(h)p_\lambda,\quad (x^-_{i})^{\lambda(h_i)+1}p_\lambda =0,\ \ i\in I,\ h\in\lie h. \end{equation} 
    \end{prop}
  \begin{pf} For $x\in\lie g$, we have
  $$x(u\otimes v)=[x,u]\otimes v+ u\otimes xv,\ \ u\in\bu(\lie g\otimes A), \ \ v\in V.$$ Since the
  adjoint action of $\lie g$ on $\lie g\otimes A$ (and hence on $\bu(\lie g\otimes A)$) is locally finite, it is
  immediate that $P_A(V)\in\cal I_A$. The proof that it is projective is  standard.
  It is clear that the element $p_\lambda\in P_A(V(\lambda))$ satisfies the relations in \eqref{defproj}
   and the fact that they are the defining relations follows by using the isomorphism in \eqref{isov}.\end{pf}
For $\nu\in P^+$ and $V\in\Ob\cal I_A$, let $V^\nu\in\Ob\cal I_A$ be the unique  maximal $\lie g\otimes A$ quotient of
$V$ satisfying \begin{equation}\label{defglobweylext}\wt(V^\nu)\subset \nu-Q^+,\end{equation} or equivalently, $$V^\nu= V/\sum_{\mu\nleq\nu}\bu(\lie g\otimes A)V_\mu.$$

A morphism $\pi:V\to V'$ of objects in $\cal I_A$ clearly induces a morphism $\pi^\nu: V^\nu\to (V')^\nu$. Let $\cal I_A^\nu$ be the full subcategory of objects $V\in\cal I_A$ such that $V=V^\nu$. It follows from the theory of finite--dimensional representations of  simple Lie algebras that \begin{equation}\label{finiteset} V\in\cal I^\nu_A\implies \#\wt V<\infty.\end{equation}

The following is immediate.
\begin{cor} Let $\nu\in P^+$ and $V\in\cal I_A^\nu$. Then $P_A(V)^\nu$ is a projective object of $\cal I_A^\nu$.
\end{cor}

\subsection{}
For $\lambda\in P^+$, set $$W_A(\lambda)=P_A(V(\lambda))^\lambda,$$
 and let $w_\lambda$ be the image of $p_\lambda$ in $W_A(\lambda)$.
The following proposition is essentially an immediate consequence of Proposition \ref{projelem} and gives
 an alternative definition of $W_A(\lambda)$ via generators and relations.  In fact this was the original
  definition given in \cite{CPweyl} when $A$ is the ring of Laurent polynomials and later generalized in \cite{FL}.

\begin{prop}\label{elemweyl} For  $\lambda\in P^+$, the module $W_A(\lambda)$ is  generated by  $w_\lambda$ with defining relations: \begin{equation}\label{weyldef}
(\lie n^+\otimes A ) w_\lambda=0,\quad hw_\lambda=\lambda(h)w_\lambda,\quad (x^-_{i})^{\lambda(h_i)+1}w_\lambda =0,\ \ i\in I,\  h\in\lie h.
\end{equation} 
\end{prop}
\begin{pf} Since $\wt W_A(\lambda)\subset\lambda-Q^+$ it follows that $(\lie n^+\otimes A) w_\lambda=0$. The other relations are clear since they are already satisfied by $p_\lambda$. To see that these  are all the relations, let $W'_A(\lambda)$ be the module generated by an element $w_\lambda$ with the relations in \eqref{weyldef}. By Proposition \ref{projelem} we see that $W'_A(\lambda)$ is a quotient of $P_A(V(\lambda))$. On the other hand $\wt(W'_A(\lambda))\subset\lambda-Q^+$ which implies that $W'_A(\lambda)$ satisfies \eqref{defglobweylext}. It follows by the maximality of $W_A(\lambda)$ that $W'_A(\lambda)$ is a quotient of $W_A(\lambda)$ and the proposition is  proved.
\end{pf}

Set \begin{gather*}\ann_{\lie g\otimes A}(w_\lambda)=\{u\in\bu(\lie g\otimes A): uw_\lambda=0\},\  \
 \ann_{\lie h\otimes A}(w_\lambda) =\ann_{\lie g\otimes A}(w_\lambda)\cap\bu(\lie h\otimes A).\end{gather*}
 Clearly $\ann_{\lie h\otimes A}(w_\lambda)$ is an ideal in $\bu(\lie h\otimes A)$ and we
 denote by
$\ba_\lambda$ the quotient of $\bu(\lie h\otimes A)$ by the ideal  $\ann_{\lie h\otimes A}(w_\lambda)$.

\subsection{} 
Regard $W_A(\lambda)$ as a right module for $\lie h\otimes A$ as follows:
 $$(uw_\lambda)(h\otimes a)=u(h\otimes a)w_\lambda,\ \ u\in\bu(\lie g\otimes A),\ h\in\lie h, a\in A.$$
 To see that this map is well defined, one must prove that:
\begin{gather*}(\lie n^+\otimes A)(h\otimes a)w_\lambda=0,\ \
(h'-\lambda(h'))(h\otimes a)w_\lambda=0,\\
(x_i^-)^{\lambda(h_i)+1}(h\otimes a)w_\lambda=0,\end{gather*} for all $i\in I$, $a\in A$ and $h,h'\in\lie h$. The first two are obvious. The
third follows from the fact that $x_i^+((h\otimes a)\otimes v_\lambda)=0$ and that  $W_A(\lambda)\in\cal I_A$. Thus, we have proved that
$W_A(\lambda)$ is a bi--module for the pair $(\lie g\otimes A, \lie h\otimes A)$.

For all $\mu\in P$, the subspaces $W_A(\lambda)_\mu$ are $\lie
h\otimes A$--submodules for both the  left and right actions and $$\ann_{\lie h\otimes A}(w_\lambda)=\{u\in\bu(\lie h\otimes A): w_\lambda
u=0=uw_\lambda\}=\{u\in\bu(\lie h\otimes A): W_A(\lambda)u=0\}.$$ Then $W_A(\lambda)$ is a $(\lie g\otimes A,\ba_\lambda)$--bimodule and each subspace $W_A(\lambda)_\mu$ is a right
$\ba_\lambda$--module. Moreover $W_A(\lambda)_\lambda$ is a $\ba_\lambda$--bimodule and we have an isomorphism of bimodules, $$W_A(\lambda)_\lambda\cong\ba_\lambda.$$

Let $\mode\ba_\lambda$ be the category of left $\ba_\lambda$--modules. Let $\bw^\lambda_A:\mode\ba_\lambda\to I_A^\lambda$ be  the right exact functor  given by
 $$\bw_A^\lambda M=W_A(\lambda)\otimes_{\ba_\lambda} M,\ \qquad \ \bw_A^\lambda f=1\otimes f,$$ where $M\in\mode\ba_\lambda$ and $f\in\Hom_{\ba_\lambda}(M,M')$ for some $M'\in\mode\ba_\lambda$.
Note that since $W_A(\lambda)\in \cal I_A$, it is clear that the $\lie g$--action on $\bw_A^\lambda M$ is also locally finite and so
$\bw_A^\lambda M\in\Ob\cal I_A^\lambda$. The preceding discussion also shows that  $$\bw^\lambda_A\ba_\lambda\cong_{\lie g \otimes A} W_A(\lambda),\qquad \ (\bw_A^\lambda
M)_\mu \cong W_A(\lambda)_\mu\otimes_{\ba_\lambda} M,\ \ \mu\in P,\ \ M\in\mode\ba_\lambda.
 $$

\subsection{} 
 \begin{lem}\label{universal} For all $\lambda\in P^+$ and $V\in\cal I_A^\lambda$ we have  $\ann_{\lie h\otimes A}(w_\lambda)V_\lambda=0$.
  \end{lem}
  \begin{pf} By Lemma \ref{superelem} and Proposition \ref{elemweyl} we see
  that  given $v\in V_\lambda$ there exists a morphism of $\lie g\otimes A$--modules
 $W_A(\lambda)\to \bu(\lie g\otimes A)v$ which maps $w_\lambda\to v$. Hence
 $uv=0$ for all $u\in\ann_{\bu(\lie h\otimes A)}(w_\lambda)$
  \end{pf}
  As a consequence of the Lemma we see that
   the left action of $\bu(\lie h\otimes A)$ on $V_\lambda$ induces a left
   action of $\ba_\lambda$ on $V_\lambda$ and we denote  the resulting
    $\ba_\lambda$--module by $\br^\lambda_A V$. Given
     $\pi\in\Hom_{\cal I_A^\lambda}(V,V')$ the restriction of $\pi_\lambda:V_\lambda\to V'_\lambda$ is a morphism of $\ba_\lambda$--modules and
    $$V\to\br^\lambda_A V,\ \ \pi\to \br^\lambda_A \pi=\pi_\lambda $$ defines a   functor $\br^\lambda_A:\cal I_A^\lambda\to\mode\ba_\lambda$ which is exact since restriction $\pi$ to a weight space is exact.
  If $M\in\Ob\mode\ba_\lambda$, we have an isomorphism of left $\ba_\lambda$--modules, $$\br^\lambda_A\bw_A^\lambda M=(\bw_A^\lambda M)_\lambda=W_A(\lambda)_\lambda\otimes_{\ba_\lambda}M\cong w_\lambda\ba_\lambda\otimes_{\ba_\lambda} M\cong M,$$
   and hence an isomorphism of functors $\id_{\ba_\lambda}\cong \br^\lambda_A\bw_A^\lambda$.

\subsection{}
 \begin{prop}\label{defeta} Let $\lambda\in P^+$ and  $V\in\cal I_A^\lambda$. There exists a canonical map of $\lie g\otimes A$--modules $\eta_V: \bw_A^\lambda\br^\lambda_A V\to V$ such that $\eta:\bw_A^\lambda\br^\lambda_A\Rightarrow\id_{\cal I_A^\lambda}$ is a natural transformation of functors and $\br^\lambda_A$ is a right adjoint to $\bw_A^\lambda$. \end{prop}
\begin{pf}   Regard $W_A(\lambda)\otimes_\bc V_\lambda$
 as a left $\lie g\otimes A$--module via the action of
 $\lie g\otimes A$ on $W_A(\lambda)$.
Lemma \ref{superelem} implies that
the assignment  $W_A(\lambda)\otimes_\bc V_\lambda\to V$ given by
  $gw_\lambda\otimes v\to gv$ is a well--defined map of left $\lie g\otimes A$--modules.
  To see that this map factors through to  a map $\eta_V:\bw_A^\lambda V_\lambda\to V$ it suffices to observe that
   $$gw_\lambda(h\otimes a)\otimes v- gw_\lambda\otimes(h\otimes a)v= g(h\otimes a)w_\lambda\otimes v-gw_\lambda\otimes (h\otimes a)v\mapsto 0 $$
   for all  $g\in\bu(\lie g\otimes A)$, $h\in\lie h$ and $a\in A$.  It is now clear that the collection $\{\eta_V; V\in\Ob\cal I_A^\lambda\}$ defines a natural transformation $\eta:\bw_A^\lambda\br^\lambda_A\Rightarrow\id_{\cal I_A^\lambda}$.

 To check that $\Res$ is right adjoint to $\bw_A^\lambda$ we must prove that there exists a natural isomorphism of abelian groups  $$\tau=\tau_{M,V}:\Hom_{\cal I_A^\lambda}(\bw_A^\lambda M, V)\cong\Hom_{\ba_\lambda}(M,\br^\lambda_A V ),$$ for all $M\in\mode\ba_\lambda$ and $V\in\cal I_A^\lambda$, such that the the following diagram commutes for all  $f\in\Hom_{\ba_\lambda}(M,M')$, $\pi\in\Hom_{\cal I_A^\lambda}(V,V')$:
 \vskip 24pt

$\begin{CD}\Hom_{\cal I_A^\lambda}(\bw_A^\lambda M',V)@>\bw_A^\lambda f^*>>\Hom_{\cal I_A^\lambda}(\bw_A^\lambda M,V)@>\pi_*>>\Hom_{\cal
I_A^\lambda}(\bw_A^\lambda M , V')\\ @VV\tau V @VV\tau V @VV\tau V \\  \Hom_{\ba_\lambda}( M',\br^\lambda_A
V)@>f^*>>\Hom_{\ba_\lambda}(M,\br^\lambda_A V)@>\Res\pi_*>>\Hom_{\ba_\lambda}(M,\br^\lambda_A V').\end{CD}$

\vskip 24pt

Define $\tau_{M,V}$ by $$\tau_{M,V}(\pi)=\pi_\lambda.$$ Since $\bw_A^\lambda M$ is generated by $M$ as a $\lie g\otimes A$--module, it follows
that $\tau(\pi)=\tau(\pi')$ implies $\pi=\pi'$. For $f\in\Hom_{\ba_\lambda}(M,\br^\lambda_A V)$ it is easily seen that
$$\tau_{M,V}(\eta_V\circ\bw_A^\lambda f)=f,$$ and hence $\tau$ is an isomorphism. The fact that the diagram commutes is straightforward.
  \end{pf}
The following is a standard consequence of properties of  adjoint functors.
  \begin{cor} The functor $\bw_A^\lambda$ maps projective objects to projective objects.
  \end{cor}

\subsection{} 
The next result gives a categorical definition of $\bw_A^\lambda M$.

\begin{thm}\label{cat} Let $V\in\cal I_A^\lambda$. Then $V\cong\bw_A^\lambda\br^\lambda_A V$ iff for all $U\in\cal I_A^\lambda$ with $U_\lambda=0$, we have \begin{equation}\label{altweyldef} \Hom_{\cal I_A^\lambda}(V,U)=0,\ \ \Ext^1_{\cal I_A^\lambda}(V, U)=0.\end{equation}
\end{thm}
\begin{pf} Suppose first that $M\in\mode\ba_\lambda$. Then $(\bw^\lambda_AM)_\lambda=w_\lambda\otimes M$ generates $\bw_A^\lambda M$ and hence
$$\Hom_{\cal I_A^\lambda}(\bw_A^\lambda M, U)=0, \ \  {\rm{if}}\ \  U_\lambda=0.$$
Let  $$P_1\to P_0\to M\to 0$$ be a right  exact sequence of modules in $\mode\ba_\lambda$, with
$P_0,P_1$   projective  and consider  the corresponding right exact sequence
$$\bw_A^\lambda P_1\to\bw_A^\lambda P_0\to \bw_A^\lambda M\to 0$$ in $\cal I_A^\lambda$. Let  $K$ be  the image of $\bw_A^\lambda P_1$
in $\bw_A^\lambda P_0$  (or equivalently the kernel of  $\bw_A^\lambda P_0\to \bw^\lambda_A M$). Then $K$ is generated as $\bu(\lie g\otimes
A)$--module by $K_\lambda$  and hence $\Hom_{\cal I_A^\lambda}(K, U)=0$ if $U\in\cal I_A^\lambda$ and $U_\lambda =0$.  By Corollary \ref{defeta} we see that $\bw_A^\lambda
P_0$ is projective and it now follows by applying $\Hom_{\cal
I_A^\lambda}(-, U)$ to the short exact sequence $$0\to K\to \bw_A^\lambda P_0\to\bw_A^\lambda M\to 0.$$ that $\Ext^1_{\cal I_A^\lambda}(\bw_A^\lambda M, U)=0$.

Conversely suppose that we are given $V\in\cal I_A^\lambda$ satisfying \eqref{altweyldef}.
 Let $V'=\bu(\lie g\otimes A)V_\lambda$ and note that $$V/V'\in\cal I_A^\lambda,\ \ (V/V')_\lambda=0.$$
 It follows from \eqref{altweyldef} that   $$\Hom_{\cal I_A^\lambda}(V,V/V')=0.$$
  This proves that $V=V'=\bu(\lie g\otimes A)V_\lambda$ and hence that the map $\eta_V: \bw_A^\lambda\br_A^\lambda V\to V$ defined in Proposition~\ref{defeta}
   is surjective. Moreover if we set $U=\ker\eta_V$, then we have $\br^\lambda_A U=0$.
 Consider the   short exact sequence $$0\to U\to\bw_A^\lambda V_\lambda\to V\to 0.$$
 Applying $\Hom_{\cal I_A^\lambda}(-, U)$  now gives $$0\to\Hom_{\cal I_A^\lambda}(U,U)\to 0,$$ and hence $U=0$ and the proof is complete.
\end{pf}

\begin{cor} The functor $\bw_A^\lambda$ is exact iff for all $U\in\cal I_A^\lambda$ with $U_\lambda=0$, we have \begin{equation}\label{equivexact} \Ext^2_{\cal I_A^\lambda}(\bw_A^\lambda M, U)=0\ \ \forall \ M\in\mode\ba_\lambda.\end{equation}
\end{cor}
\begin{pf} Assume that \eqref{equivexact} is satisfied. Let $0\to M''\to M\to M'\to 0$ be a short exact sequence of modules in $\mode\ba_\lambda$ and consider the induced  short exact sequence $$0\to K\to \bw_A^\lambda M\to \bw_A^\lambda M'\to 0.$$
Apply $\Hom(-, U)$ to the preceding short exact sequence and using Theorem~\ref{cat} and \eqref{equivexact} we find that $$\Hom_{\cal
I_A^\lambda}(K,U)=0,\ \ \Ext^1_{\cal I_A^\lambda}(K, U)=0,\ \ \forall\ \ U\in\Ob\cal I_A^\lambda \text{ with } U_\lambda=0$$ Hence
$K\cong\bw_A^\lambda K_\lambda$. Applying the functor $\br^\lambda_A$ and using the fact that $\br^\lambda_A\bw_A^\lambda$ is naturally
isomorphic to the identity functor,  we see that if $V$ is  the kernel of $\bw_A^\lambda M''\to K$ then $V_\lambda=0$. Applying $\Hom_{\cal
I_A^\lambda}(- ,V)$ to the short exact sequence $$0\to V\to \bw_\lambda M''\to K\to 0,$$ proves that $V=0$.

For the converse, suppose that $\bw_A^\lambda$ is exact. Let $M\in\Ob\mode\ba_\lambda$ and let $P\in\Ob\mode\ba_\lambda$ be projective  such that we have an exact sequence $0\to M'\to P\to M\to 0.$ This gives us $$0\to\bw_A^\lambda M'\to\bw_A^\lambda P\to\bw_A^\lambda M\to 0.$$ Applying $\Hom_{\cal I_A^\lambda} (-, U)$ with $U\in\cal I^\lambda_A$, $U_\lambda=0$ and recalling that $\bw_A^\lambda P$ is projective in $\cal I_A^\lambda$ we get a piece of the long exact sequence  $$0\to \Ext^2(\bw_A^\lambda M, U)\to 0,$$ and the converse is established.
\end{pf}

\section{The structure of  $W_A(\lambda)$}

\subsection{}
We begin by proving that the construction of  $W_A(\lambda)$ is functorial in $A$. Assume that $B$ is a commutative associative algebra and let $f: A\to B$ be a homomorphism of algebras. Then $(1\otimes f):\lie g\otimes A\to\lie g\otimes B$ is a homomorphism of Lie algebras and given any $\lie g\otimes B$--module $V$ we can regard it as a $\lie g\otimes A$--module via $f$ and we denote this module by $f^*V$.
 \begin{prop}\label{canonical} Let $\lambda\in P^+$ and let $f:A\to B$ be a homomorphism of associative algebras. Then $f$ induces a canonical homomorphism $f_\lambda:\ba_\lambda\to \bb_\lambda$ of associative algebras and   a canonical map of ($\lie g\otimes A, \ba_\lambda)$-bimodules $f_\lambda^*:W_A(\lambda)\to f^*(W_B(\lambda))$. Moreover, $f_\lambda$ and $f_\lambda^*$ are surjective if $f$ is surjective.
\end{prop}
\begin{pf} The action of $\lie g\otimes A$ on $f^*(W_B(\lambda))$ is given by $$(x\otimes a)\circ w_{\lambda, B}=(x\otimes f(a))w_{\lambda, B}$$ and it follows immediately from Proposition \ref{elemweyl} that
 there is a well--defined map of left $\lie g\otimes A$--modules $$W_A(\lambda)\to  f^*(W_B(\lambda)),\ \  \ \ w_{\lambda, A}\to w_{\lambda, B}.$$ Since $(1\otimes f)$ maps $\lie h\otimes A$ to $\lie h\otimes B$ this is also a map of right $\bu(\lie h\otimes A)$--modules. The proof of the proposition is complete if we prove that  $$u\in\ann_{\lie h\otimes A}(w_{\lambda, A})\implies (1\otimes f)(u)\in\ann_{\lie h\otimes B}(w_{\lambda, B}).$$
 But this is clear since $$w_{\lambda,A}u=uw_{\lambda,A}\to(1\otimes f)(u)w_{\lambda, B}=w_{\lambda,B}(1\otimes f)(u).$$
\end{pf}
 Let $A,B$ and $f:A\to B$ be as in the proposition and given $M\in\mode\bb_\lambda$, let $f_\lambda^*M\in\mode\ba_\lambda$ be the corresponding $\ba_\lambda$--module.
\begin{cor}  There exists a natural morphism of $\lie g\otimes A$--modules $\bw^\lambda_Af_\lambda^*M\to f^*\bw^\lambda_BM$ which is surjective if $f$ is surjective.  In particular  we have a morphism  of $\lie g\otimes A$--modules \begin{equation}\label{pullbackiso}\bw^\lambda_Af_\lambda^*\bb_\lambda\to f^*\bw^\lambda_B\bb_\lambda\cong f^*(W_B(\lambda)),\end{equation} which is  surjective if $f$ is surjective.
 \end{cor}
 \begin{pf} It is clear that there exists a map $f^*\otimes f_\lambda^*$ of $\lie g\otimes A$--modules
 $$W_A(\lambda)\otimes_{\ba_\lambda} f_\lambda^*M= \bw^\lambda_Af_\lambda^*M\longrightarrow f^*W_B(\lambda)\otimes_{\ba_\lambda} f_\lambda^*M.$$ Composing with the map of $\lie g\otimes A$--modules, $$f^*W_B(\lambda)\otimes_{\ba_\lambda} f_\lambda^*M \to f^*\bw^\lambda_BM =f^*(W_B(\lambda)\otimes_{\bb_\lambda} M),\ \qquad u\otimes m\to u\otimes m$$ proves the corollary.
 \end{pf}

\subsection{} 
The next proposition begins an analysis of the behaviour of the modules $W_A(\lambda)$ and the functors $\bw^\lambda_A$ under tensor products. We shall assume from now on that an unadorned $\otimes$ denotes the tensor product of vector spaces over $\bc$.

\begin{prop} \label{cantensor} Let $\lambda,\mu\in P^+$. \begin{enumerit} \item There exists a homomorphism of $\lie g\otimes A$--modules   $$\tau_{\lambda,\mu}: W_A(\lambda+\mu)\to W_A(\lambda)\otimes W_A(\mu),$$ such that
$\tau_{\lambda,\mu}(w_{\lambda+\mu})=\ w_\lambda\otimes w_\mu.$
 \item The homomorphism $\Delta: \bu(\lie h\otimes A)\to\bu(\lie h\otimes A)\otimes \bu(\lie h\otimes A)$ induces a canonical homomorphism $\Delta_{\lambda,\mu}:\ba_{\lambda+\mu}\to\ba_\lambda\otimes\ba_\mu$ and $$\Delta_{\lambda,\mu}=\sigma_{\mu, \lambda} \circ \Delta_{\mu, \lambda},\ \ (1\otimes \Delta_{\mu,\nu})\circ \Delta_{\lambda,\mu+\nu}=( \Delta_{\lambda,\mu} \otimes 1)\circ \Delta_{\lambda+\mu,\nu},\ \ \nu\in P^+.$$
 where $\sigma_{\lambda,\mu}: \ba_{\lambda} \otimes \ba_{\mu} \longrightarrow \ba_{\mu} \otimes \ba_{\lambda}$ denotes the flip map.
     \item The tensor product $W_A(\lambda)\otimes W_A(\mu)$ is canonically a $(\lie g\otimes A, \ba_\lambda\otimes\ba_\mu)$--bimodule  and hence also a $(\lie g\otimes A, \ba_{\lambda+\mu})$--bimodule.
         \item The map $\tau_{\lambda,\mu}$ is a map of $(\lie g\otimes A, \ba_{\lambda+\mu})$--bimodules and for $M\in\mode\ba_\lambda$,
         $N\in\mode\ba_\mu$ we have an induced map of $\lie g\otimes A$-modules
          $$\tau_{M,N}: \bw_A^{\lambda+\mu}\Delta_{\lambda,\mu}^*(M\otimes N)\to \bw_A^\lambda M\otimes \bw^\mu_A N.$$
     \end{enumerit}
 \end{prop}
 \begin{pf}  Part (i) is immediate from Proposition \ref{elemweyl}. It follows that  $$u\in\ann_{\lie h\otimes A}(w_{\lambda+\mu})\implies \Delta(u)(w_\lambda\otimes w_\mu)=0,$$ i.e., that $$\Delta(u)\in\ann_{\lie h\otimes A}(w_\lambda)\otimes \bu(\lie h\otimes A)+\bu(\lie h\otimes A)\otimes \ann_{\lie h\otimes A}(w_\mu),$$
  and hence we have an induced map $\Delta_{\lambda,\mu}:\ba_{\lambda+\mu}\to\ba_\lambda\otimes\ba_\mu$. The remaining statements in (ii)  follow from the co-commutativity and co-associativity of $\Delta$. The right action of $\ba_\lambda$ on $W_A(\lambda)$ and of $\ba_\mu$ on $W_A(\mu)$ defines  a right action of $\ba_\lambda\otimes \ba_\mu$ on $W_A(\lambda)\otimes W_A(\mu)$ in the obvious pointwise way and part (iii) now follows easily.  To prove (iv), note that we clearly have a map $$\bw^{\lambda+\mu}_A\Delta_{\lambda,\mu}^*(M\otimes N)\to \left(W_A(\lambda)\otimes W_A(\mu)\right)\otimes_{\ba_{\lambda+\mu}}\Delta_{\lambda,\mu}^*(M\otimes N).$$
  Since there exist  canonical maps of $\lie g\otimes A$--modules \begin{gather*}\left( W_A(\lambda)\otimes W_A(\mu)\right)\otimes_{\ba_{\lambda+\mu}}\Delta_{\lambda,\mu}^*(M\otimes N)\to \left( W_A(\lambda)\otimes W_A(\mu)\right)\otimes_{\ba_{\lambda}\otimes\ba_{\mu}}(M\otimes N)\end{gather*} and a map \begin{gather*}\left( W_A(\lambda)\otimes W_A(\mu)\right)\otimes_{\ba_{\lambda}\otimes\ba_{\mu}}(M\otimes N) \to \bw^\lambda_A M\otimes\bw^\mu_A N,\\ (w\otimes w')\otimes (m\otimes n)\to (w\otimes m)\otimes (w'\otimes n),\end{gather*} the result follows.
 \end{pf}

\subsection{}
 Given two commutative associative algebras $A$ and $B$ the direct sum  $C=A\oplus B$ is
canonically an associative algebra and let $p_{A}$ (resp. $p_B$)  be the projection onto $A$ (resp. $B$). By Proposition \ref{canonical}  any
$M\in\mode\ba_\lambda$ (resp. $N\in\mode\bb_\mu$)  can be regarded as a module for $\bc_\lambda$ (resp. $\bc_\mu$) and hence the tensor
product $M\otimes N$ can be viewed as a module for $\bc_\lambda\otimes\bc_\mu$. Pulling this module back by $\Delta_{\lambda,\mu}$ we get
 a $\bc_{\lambda+\mu}$--module
  which by abuse of notation, we shall just denote by $M\otimes N$
 and  we shall see that the context is such that no confusion arises from this abuse of notation.
The following is immediate from Corollary \ref{canonical} and  Proposition \ref{cantensor}(iv).
 \begin{cor}\label{coro-3-3} For $M \in \mode \ba_{\lambda}, N \in \mode\bb_{\mu}$, there exists a surjective homomorphism of $\lie g\otimes C$--modules $$\bw_C^{\lambda+\mu} (M \otimes N)\twoheadrightarrow \bw_A^{\lambda}M\otimes \bw_{B}^{\mu}N.$$
 \end{cor}

\subsection{}
\begin{thm} \label{fingen} Assume that  $A$ is a finitely generated  algebra.
\begin{enumerit} \item[(i)]For all  $\lambda\in P^+$, the algebra $\ba_\lambda$ is finitely generated and $W_A(\lambda)$ is a finitely generated right $\ba_\lambda$--module.
\item[(ii)] If  $M\in\mode\ba_\lambda$ is a  finitely generated (resp. finite--dimensional) then $\bw_A^\lambda M$ is a finitely generated (resp. finite-dimensional) $\lie g\otimes A$--module.
\item[(iii)] Suppose that $A$ and $B$ are finite--dimensional commutative, associative algebras and let $\lambda,\mu\in P^+$. For
$M\in\mode\ba_\lambda$, $N\in\mode \bb_\mu$ with $\dim M<\infty$ and $\dim N<\infty$ we have,
  $$\bw_{A\oplus B}^{\lambda+\mu}(M\otimes N)\cong \bw^\lambda_A M\otimes \bw^\mu_B N,$$ as $\lie g\otimes (A\oplus B)$--modules.
    \end{enumerit}
\end{thm}
We prove the theorem in the rest of the section.

\subsection{} 
Let $u$ be an indeterminate and for  $a\in A$, $\alpha\in R^+$, define a power series $\bop_{a,\alpha}(u)$ in  $u$ with coefficients in $\bu(\lie h\otimes
A)$ by $$\bop_{a,\alpha}(u)=\exp\left(-\sum_{r=1}^\infty \frac{h_\alpha\otimes a^r}{r}u^r\right).$$ For $s\in\bz_+$, let $p_{a,\alpha}^s$ be the
coefficient of $u^s$ in $\bop_{a,\alpha}(u)$.
The following formula is proved in \cite{G} in the case when $A$ is the  polynomial ring $\bc[t]$ and $a=t$. Applying the Lie algebra
homomorphism  $$\lie g\otimes \bc[t]\to \lie g\otimes A,\ \ \ x\otimes t^r\to x\otimes a^r,\ \ r\in\bz_+,\ \ x\in\lie g,$$ gives the result for $\lie g\otimes ^{}A$.
\begin{lem} \label{garland} Let $r\in\bz_+$. Then, \begin{gather*}(x^+_\alpha\otimes a)^r(x^-_\alpha\otimes 1)^{r+1}-\sum_{s=0}^r (x_\alpha^-\otimes a^{r-s})p^{s}_{a,\alpha}
\in\bu(\lie g\otimes A)(\lie n^+\otimes A),\\ (x^+_\alpha\otimes a)^{r+1}(x^-_\alpha\otimes 1)^{r+1}-p^{r+1}_{a,\alpha}
\in\bu(\lie g\otimes A)(\lie n^+\otimes A)
\end{gather*}\hfill\qedsymbol
\end{lem}

\subsection{} 
Part (i) of the theorem was proved in the case when $A$ is the polynomial ring in one variable in \cite{CPweyl}. The proof in the general case is very similar, and we only give a brief sketch here.  Let $a_1,\cdots ,a_m$ be a set of generators for $A$. Using the defining relations of
$W_A(\lambda)$ and Lemma \ref{garland}, we see that $$(x^+_i\otimes a_k)^{n_i}(x^-_i\otimes 1)^{n_i+1}w_\lambda=\sum_{s=0}^{n_i} (x_i^-\otimes
a_k^{n_i-s})p^{s}_{a_k,\alpha_i}w_\lambda=0$$ for all $i\in I$, $1\le k\le m$ and $n_i=\lambda(h_i). $
Applying $x_i^+\otimes a$, $a\in A$,  to both sides of the equation, we  get  $$\left(h_i\otimes aa_k^{n_i}+\sum_{s=1}^{n_i} (h_i\otimes aa_k^{n_i-s})p^s_{a_k,\alpha_i}\right)w_\lambda=0.$$ It is now straightforward to see by using an iteration of this argument  that for all $i\in I$, $(r_1,\cdots, r_m)\in\bz_+^m$, we have $$h_i\otimes (a_1^{r_1}\cdots a_m^{r_m})w_\lambda=H(i,r_1,\cdots, r_m)w_\lambda$$ for some $H(i,r_1,\cdots ,r_m)$ in the  subalgebra of $\bu(\lie h\otimes A)$ generated by the elements of the set $$  \{h_i\otimes a_1^{s_1}\cdots a_m^{s_m}: 0\le s_\ell\le n_i,\ 1\le \ell\le m, \ i\in I\}.$$ In other words, we have proved that $\ba_\lambda$ is the quotient of a finitely generated algebra.

Let $\{\beta_1,\cdots,\beta_N\}$ be an enumeration of $R^+$ and set $$S=\{a_1^{s_1}\cdots a_m^{s_m}:(s_1,\cdots,s_m)\in\bz_+^M\}.$$
    Using the PBW theorem, we see that  elements of the set,  \begin{equation}\label{gen} \left\{(x^-_{\beta_{i_1}}\otimes b_1)\cdots (x^-_{\beta_{i_\ell}}\otimes b_\ell)w_\lambda: 1\le i_1\le \cdots\le i_\ell\le N,\ \ \ell\in\bz_+,\  b_1,\cdots,b_\ell\in S\right\}\end{equation}  generate $W_A(\lambda)$ as a right module for $\ba_\lambda$.  Using Lemma \ref{garland} and the defining relations for $W_A(\lambda)$ we see that
 $$(x^+_\alpha\otimes a_r)^{n_\alpha}(x^-_\alpha\otimes 1)^{n_\alpha+1}w_\lambda=\sum_{s=0}^{n_\alpha} x_\alpha^-\otimes a_r^{n_\alpha-s}p^s_{a_r,\alpha}w_\lambda=0,\ \ 1\le r\le m,$$ for all $\alpha\in R^+$ and $n_\alpha=\lambda(h_\alpha)$. That implies  $$(x^-_\alpha\otimes a_r^s)w_\lambda\in{\rm{sp}} \{(x^-_\alpha\otimes a_r^\ell)w_\lambda\ba_\lambda: \ \ 0\le \ell<\lambda(h_\alpha)\}.$$
 Applying $h_\alpha\otimes a^k_p$ with $r\ne p$ to the preceding equation gives,
\begin{gather*}(x^-_\alpha\otimes a_r^sa_p^k)w_\lambda \in{\rm{sp}}\{(x^-_\alpha\otimes a_r^\ell
a_p^{k})w_\lambda\ba_\lambda: \ \ 0\le \ell<\lambda(h_\alpha)\}\\
\subset{\rm{sp}} \{x^-_\alpha\otimes a_r^\ell a_p^{\ell'} W_A(\lambda)_\lambda,\ \ 0\le \ell,\ell'<n_\alpha\}.\end{gather*} It is now clear that more
generally we have $$(x^-_\alpha\otimes A)w_\lambda\subset {\rm{sp}} \{(x^-_\alpha\otimes (a_1^{r_1}\cdots a_m^{r_m})w_\lambda\ba_\lambda: 0\le
r_\ell<n_\alpha)\}.$$ An induction on the length of the monomials in \eqref{gen} identical to the one used in \cite{CPweyl} now proves that $W_A(\lambda)$ is a finitely generated $\ba_\lambda$--module.  Part (ii) of the theorem is now immediate  by using \eqref{finiteset}.

\subsection{}
To prove (iii), we begin with the following refinement of Theorem \ref{cat}.
\begin{prop}
\label{refi-defn}
 \begin{enumerit}
 \item[(i)] Let $\lambda,\nu\in P^+$ be such that $\lambda\nleq\nu$ and $\nu\nleq\lambda$. Let $U\in\cal I^\nu_A$ be irreducible and assume that   $U_\nu\ne 0$. Then $$\Ext^m_{\cal I_A}(\bw^\lambda_AM, U)=0,\ \ m=0,1,$$ for all
$M\in\Ob\mode\ba_\lambda$.
\item[(ii)] Let  $V\in\cal I^\lambda_A$ be such that $\dim V_\lambda<\infty$. Then $\bw^\lambda_A\br^\lambda_AV_\lambda\cong V$ iff \begin{equation}\label{altweyldefa} \Ext^m_{\lie g\otimes A}(V, U)=0,\ \ m=0,1\end{equation}
  for all $U\in\Ob \cal I^\lambda_A$ with  $\dim U<\infty$ and  $U_\lambda=0$.
\end{enumerit}
\end{prop}
\begin{pf}  For (i), observe that since $U$ is irreducible any non--zero morphism
  $\eta:W_A(\lambda)\to U$ must be surjective. But this is impossible since
  $(\bw^\lambda_AM)_\nu=0$. Suppose next that $$0\to U\to V\to \bw^\lambda_AM\to 0$$ is a short exact sequence of objects in $\cal I_A$.
  Then $$V_\lambda\ne 0,\ \ \ \  \wt V\subset (\nu-Q^+)\cup(\lambda-Q^+),$$ and since $\lambda\nleq\nu$ we
  see that  $(\lie n^+\otimes A)V_\lambda=0.$ Set $V'=\bu(\lie g\otimes A)V_\lambda$ so
   that $\wt V\subset\lambda-Q^+$. To prove that the sequence splits, it suffices to
   prove that $$V'\cap U=\{0\}.$$ Otherwise  since $U$ is irreducible  we would have $U\cap V'=U$ which
    would imply that $\nu\in \wt V'$ contradicting $\nu\nleq\lambda$.

A simple induction on the length of $U$ shows that it suffices to  to prove that $\bw^\lambda_AV_\lambda\cong V$ if   \eqref{altweyldefa} holds for all irreducible modules $U\in\Ob\cal I^\lambda_A$ with $U_\lambda=0$. As in the proof of Theorem \ref{cat} we have $V=\bu(\lie g\otimes A)V_\lambda$ and hence a short exact sequence $$0\to K\to \bw^\lambda_A V_\lambda\to V\to 0.$$ By part (ii) of Theorem~\ref{fingen} we have $\dim\bw^\lambda_AV_\lambda<\infty$ and hence we have $$
 \dim K<\infty,\ \ K_\lambda=0.$$ If $K\ne 0$,
 then $\Hom_{\lie g\otimes A}(K,U)\ne 0$ for some irreducible module $U\in\cal I^\lambda_A$ with $U_\lambda=0$. Applying $\Hom_{\cal I^\lambda_A}(-, U)$ and using the fact that $\Hom_{\lie g\otimes A}(\bw^\lambda_A, U)=0$, we get $$0\to\Hom_{\lie g\otimes A}(K,U)\to\Ext^1_{\lie g\otimes A}(V,U)$$
 which is impossible since $V$ satisfies \eqref{altweyldefa}. Hence $K=0$ and the proof of (ii) is complete.
\end{pf}

\subsection{}
 The proof of part(iii) of the Theorem is completed as follows. By Corollary~\ref{coro-3-3} we have a surjective map of $\lie g\otimes (A\oplus B)$--modules,
$$\bw_{A\oplus B}^{\lambda+\mu}(M\otimes N)\longrightarrow \bw^\lambda_AM\otimes \bw^\mu_BN\to 0.$$  To prove that it is an isomorphism it suffices by Proposition ~\ref{refi-defn}(ii)  to prove that
$$\Ext^{m}_{\cal I^{\lambda+\mu}_{A \oplus B}}(\bw^\lambda_AM\otimes \bw^\mu_BN, U) = 0,\ \ m=0,1,$$
 for all irreducible  $U\in\Ob\cal I^{\lambda+\mu}_{A\oplus B}$ with $U_{\lambda + \mu} = 0$.
 By Proposition ~\ref{ideal-irred} we may write such a module as a tensor product,
  $$U \cong U_A \otimes U_B,\ \ U_A\in\Ob\cal I_A, \ \ U_B\in\Ob\cal I_B, $$
  where $U_A$ and  $U_B$ are irreducible. Let $\nu_A$ (resp. $\nu_B$)  be the highest weight of $U_A$ (resp. $U_B$) and note that $\nu_A+\nu_B\in\wt U\subset\lambda+\mu-Q^+$.
  Since $\bw^\lambda_AM$, $ \bw^\mu_BN$ and $U$ are all finite--dimensional modules for finite--dimensional Lie algebras, we have for $m=0,1$,
   \begin{gather*} \Ext^m_{\lie g\otimes(A\oplus B)}(\bw^\lambda_AM\otimes \bw^\mu_BN, U)
 \cong \Ext^{m}_{\cal I^{\lambda+\mu}_{A \oplus B}}(\bw^\lambda_AM\otimes \bw^\mu_BN, U),\\ \Ext^m_{\lie g\otimes A}(\bw^\lambda_AM, U_A)
 \cong \Ext^{m}_{\cal I^\lambda_{A }}(\bw^\lambda_AM,U_A),\qquad
 \Ext^m_{\lie g\otimes B}( \bw^\mu_BN, U_B)
 \cong \Ext^{m}_{\cal I^\lambda_{b }}( \bw^\mu_BN,U_B).\end{gather*}
By Proposition \ref{kunneth} it suffices to prove that either
 \begin{equation}\label{vanishing}\Ext^{m}_{\cal I^\lambda_A}(\bw^\lambda_AM, U_A) = 0,\ \\ \ {\rm{or}} \ \  \Ext^{m}_{\cal I^\mu_B}( \bw^\mu_BN, U_B) = 0,\ \ m=0,1. \end{equation} If $U_A\in\Ob\cal I^\lambda_A$ or $U_B\in\Ob\cal I^\nu_B$  then \eqref{vanishing} follows from  Proposition ~\ref{refi-defn}(ii).  Otherwise  we have $$\nu_A\nleq\lambda,\ \qquad\ \nu_B\nleq \mu.$$ Since $\nu_A+\nu_B<\lambda+\mu$, it follows now that $\lambda\nleq\nu_A$ and now \eqref{vanishing} follows from
  Proposition ~\ref{refi-defn}(i).

\section{Further results on tensor products } Throughout this section, we assume that $A$ is finitely generated.

\subsection{}
Let $\rm{irr}\mode\ba_\lambda$ be the set of irreducible representations of $\ba_\lambda$. Since $\ba_\lambda$ is a commutative finitely generated algebra it follows that if $M\in\rm{irr}\mode\ba_\lambda$ then $\dim M=1$. By Theorem \ref{fingen} we see that $$\dim\bw_A^\lambda M<\infty,\ \ \br_A^{\lambda}\bw^\lambda_AM =M,\ \text{ for }\ M\in\rm{irr}\mode\ba_\lambda,$$ and we denote by
$\bv^\lambda_AM$ the unique irreducible quotient of $\bw_A^\lambda M$ (see Lemma \ref{superelem}).
  It now follows from  Lemma \ref{ideal} and Lemma \ref{sscond}  that there exists an ideal of finite--codimension $\tilde K^\lambda_M$ of $A$ such
    that $\lie g\otimes A/\tilde K^\lambda_M$ is a semisimple Lie algebra and $$(x\otimes a)\bv_A^\lambda M=0\ \
    \forall\ \ x\in\lie g,\ \ a\in \tilde{K}^\lambda_M.$$
    Suppose that $M\in\mode\ba_\lambda$ is finite dimensional of length $r$, $M_1,\cdots, M_r$ be the irreducible constituents of $M$ and set $$\tilde
    K^\lambda_M=\prod_{s=1}^r\tilde K^\lambda_{M_s}.$$

\subsection{}
  The next result shows that any irreducible module in $\cal I_A^\lambda$ is isomorphic to $\bv^\mu_AM$ for some $\mu\in P^+$. \begin{lem}\label{uniqueirr} Let $\lambda\in P^+$ and assume that $V\in\cal I_A^\lambda$ is irreducible.  There exists $\mu\in P^+\cap \cal(\lambda- Q^+)$ such that $$\wt V\subset\mu- Q^+,\ \ \dim V_\mu=1.$$
 In particular, $V$ is the unique irreducible  quotient of $\bw_A^\mu\br_A^\mu V$ and hence $\dim V<\infty$. If $V'\in\Ob \cal I_A$ we have  $V\cong V'$ as $\lie g\otimes A$--modules iff $\br^\mu_AV\cong \br^{\mu'}_AV'$ as $\ba_\mu$--modules.\end{lem}
 \begin{pf} Since $V\in\cal I^\lambda_A$, it follows that there exists  $\mu\in\lambda-Q^+$ with  $$V_\mu\ne 0,\ \  \ \ (\lie n^+\otimes A)V_\mu=0.$$ It is immediate from Proposition \ref{defeta} that $V$ is a quotient of $\bw_A^\mu \br_A^\mu$.
  If  $V_\mu'=\bu(\lie h\otimes A)V_\mu$ is a proper $\lie h\otimes A$--submodule of $V_\mu$, then $V'=\bu(\lie g\otimes A)V_\mu'$ is a proper submodule of $V$ which is a contradicton. Hence $ \br_A^\mu V$ is an irreducible  $\ba_\mu$--module which implies that $\dim V_\mu=1$. Theorem \ref{fingen} now implies that $\dim \bw_A^\mu \br_A^\mu V<\infty$ and hence $\dim V<\infty$. The proof that $V$ is the unique irreducible quotient of $\bw^A_\mu \br_A^\mu V$ is standard since $\br_A^{\mu}\bw^\mu_A \br_A^\mu V\cong V_{\mu}$.
  The final statement of the lemma is  now trivial.
 \end{pf}

\subsection{}
The main result of this section is the following.
\begin{thm} \label{tensprod}Let $\lambda,\mu\in P^+$ and let $M,N$ be irreducible modules for
 $\ba_\lambda$ and $\ba_\mu$ respectively and assume that \begin{equation}\label{ass} A/\tilde K^\lambda_M\tilde K^\lambda_N\cong
A/\tilde K^\lambda_M\oplus A/\tilde K^\lambda_N.\end{equation}
 Then \begin{gather}\label{tensirr}\bv_A^{\lambda+\mu}(M\otimes N)\cong_{\lie g\otimes A} \bv_A^\lambda M\otimes \bv^\mu_A N,\qquad
  \ \tilde K^{\lambda+\mu}_{M\otimes N}=\tilde{K}^\lambda_M\tilde{K}^\mu_N, \\ \label{tensweyl}
\bw_A^{\lambda+\mu}(M\otimes N)\cong_{\lie g\otimes A} \bw_A^\lambda M\otimes \bw^\mu_A N.\end{gather}
\end{thm}

\subsection{}
 To prove \eqref{tensirr} recall
 that $M\otimes N$ is an irreducible $\ba_\lambda\otimes\ba_\mu$--module with the action being pointwise and hence also an irreducible $\ba_{\lambda+\mu}$--module (via
$\Delta_{\lambda,\mu}$). By Lemma \ref{uniqueirr} we see that it suffices to prove that $\bv_A^\lambda M\otimes \bv^\mu_A N$ is the irreducible
$\lie g\otimes A$ quotient of $\bw^{\lambda+\mu}_A(M\otimes N)$. Clearly,  $\bv_A^\lambda M\otimes \bv^\mu_A N$ is an irreducible module for the
semisimple Lie algebra $\lie g\otimes (A/\tilde K^\lambda_M\oplus A/\tilde K^\lambda_N)$   and hence using
  \eqref{ass} it is an irreducible module for $\lie g\otimes A/\tilde K^\lambda_M\tilde K^\lambda_N$
  and so for  $\lie g\otimes A$ as well. Since
  $$\br_A^{\lambda + \mu} (\bv_A^\lambda M\otimes \bv^\mu_A N)\cong M\otimes N,$$
  we see from Lemma \ref{universal} that  $\bv_A^\lambda M\otimes \bv^\mu_A N$ is a quotient of $\bw^{\lambda+\mu}_A(M\otimes N)$ and the first isomorphism in \eqref{tensirr} is proved. For the second, observe that
by definition if $S$ is any ideal in $A$ such that $$(\lie g\otimes S)\bv_A^{\lambda}M=0,$$ then $S\subset\tilde{K}^{\lambda}_M$
and similarly for $\tilde{K}^{\mu}_N$.  One deduces easily from  \eqref{ass}
  that $\tilde{K}^\lambda_M\tilde{K}^\lambda_N$ is the largest ideal in $A$ such that
   $$(\lie g\otimes \tilde{K}^\lambda_M\tilde{K}^\lambda_N)\bv_A^\lambda M\otimes \bv^\mu_A N=0.$$ Since
    $\tilde K^{\lambda+\mu}_{M\otimes N}$  is  maximal with the property that $$(\lie g\otimes\tilde K^{\lambda+\mu}_{M\otimes N})\bv_A^{\lambda+\mu}(M\otimes N)=0 $$ we now get that $\tilde K^{\lambda+\mu}_{M\otimes N} = \tilde{K}^\lambda_M\tilde{K}^\lambda_N$.

\subsection{}
We need several results to prove \eqref{tensweyl}.  Theorem \ref{fingen} and Lemma \ref{ideal} imply  that given $\lambda\in P^+$ and $M\in\mode\ba_\lambda$ with $\dim M<\infty$, there exists
an ideal of finite codimension $K^\lambda_M$ in $A$ which is maximal with the property that
$$(\lie g\otimes K^\lambda_M)\bw_A^\lambda M=0.$$  If $0\to M'\to M\to M''\to 0,$ is a short exact sequence of modules in $\ba_\lambda$ then since the functor $\bw^\lambda_M$ is right exact, we see that
 \begin{equation}\label{product}K^\lambda_{M'}K^\lambda_{M''}\subset K^\lambda_M\subset K^\lambda_{M''}.\end{equation}
 Let $K\subset K^\lambda_M$ be an ideal in $A$ and set   $A/K=B$. It is clear that $\bw^\lambda_AM$ is a module for $\lie g\otimes B$
and since $$\br^\lambda_B\bw^\lambda_AM =M,$$ we get  by Lemma \ref{universal} that $M$ is also a  $\bb_\lambda$--module.
\begin{lem} Let $\lambda\in P^+$ and $M\mode\ba_\lambda$ be finite--dimensional.
For all ideals $K\subset K^\lambda_M$, we have an isomorphism of $\lie g\otimes A$ (or equivalently $\lie g\otimes A/K$) modules,
\begin{equation}\label{limitideal}\bw^\lambda_AM\cong \bw^\lambda_{A/K}M.\end{equation}
\end{lem}
 \begin{pf}
  By  Corollary \ref{canonical} and the discussion preceding the statement of the Lemma we see that we have a surjective map of $\lie g\otimes A$--modules $$\bw_A^\lambda M\to\bw^\lambda_B M\to 0,\ \ w_\lambda\otimes m\to w_\lambda\otimes m.$$
  On the other hand by Proposition~\ref{defeta} we have a map of $\lie g\otimes B$--modules $$\bw^\lambda_BM\cong \bw^\lambda_B\br^\lambda_B\bw^\lambda_AM\longrightarrow \bw^\lambda_AM,\  w_\lambda\otimes m\to w_\lambda\otimes m$$ and hence
  \eqref{limitideal} is proved.
  \end{pf}

\subsection{}
\begin{prop}\label{irrweyl} Let $\lambda\in P^+$ and $M\in\mode\ba_\lambda$ be finite--dimensional.
We have $$(\tilde K^\lambda_M)^{\lambda(h_\theta)}\subset K^\lambda_M.$$\end{prop}
 \begin{pf} It suffices by \eqref{product} to consider the case when $M$ is irreducible.  Using Lemma \ref{garland}
 we  see as in the proof of Theorem \ref{fingen} that
 $$0=(x^+_\theta\otimes a)(x^-_\theta)^{\lambda(h_{\theta})+1}(w_\lambda\otimes m)= \sum_{s=0}^{\lambda(h_\theta)}
 ( x^-_\theta\otimes a^{r-s}) p^s_{a,\theta}(w_\lambda\otimes m).$$ If $a\in \tilde K^\lambda_M$ then
 $(h\otimes a)(w_\lambda\otimes m)=0$ and since $p^s_{a,\theta}$ is in the subalgebra generated by
 the elements $\{h_\theta\otimes a^p:  p\in\bz_+, p>0\}$ with constant term zero, we see
  that $p^s_{a,\theta}(w_\lambda\otimes m)=0$ for all $s>0$. This implies that
   $$(x^-_\theta\otimes a^{\lambda(h_\theta)})(w_\lambda\otimes m)=0.$$  Since $[x^-_\theta,\lie n^-]=0$ we get
   $$(x^-_\theta\otimes a^{\lambda(h_\theta)})\bw_A^\lambda M=0.$$
    Since $\lie g$ is generated by $x^-_\theta$ as a $\lie g$--module the result follows.\\
\end{pf}

\subsection{}
By part (i) of the theorem  and Proposition \ref{irrweyl},
     we may choose $r\ge 1$ so that $$(\tilde K^\lambda_M)^r(\tilde K_N^\mu)^r=
     (\tilde K_{M\otimes N}^{\lambda+\mu})^r\subset(K_M^\lambda K^\mu_M)\cap K^{\lambda+\mu}_{M\otimes N}.$$
 Set $C=A/(\tilde K^\lambda_M\tilde K^\mu_N)^r$ and note that $$C=A/(\tilde K^\lambda_M)^r\oplus A/(\tilde K^\mu_N)^r.$$
   By Theorem \ref{fingen}(ii), we have an isomorphism of $\lie g\otimes A/C$--modules
    \begin{equation*}\bw_C^{\lambda+\mu}(M\otimes N)
    \cong\bw_{A/(\tilde K^\lambda_M)^r}^\lambda M\otimes\bw^\mu_{A/(\tilde K^\mu_N)^r} N.\end{equation*}
     Lemma \ref{limitideal} now proves that we have  isomorphisms of $\lie g\otimes A$--modules,
     $$\bw_C^{\lambda+\mu}(M\otimes N)\cong \bw^{\lambda+\mu}_A(M\otimes N),\ \
      \ \bw_{A/(\tilde K^\lambda_M)^r} M\cong\bw^\lambda_A M,\ \  \ \bw^\mu_{A/(\tilde K^\mu_N)^r} N\cong \bw^\mu_A N,$$ and \eqref{tensweyl} is proved.

\subsection{}
 The statement of \eqref{tensweyl} can be strengthened as follows by using Proposition~\ref{irrweyl}.
\begin{cor} Let $M\in\mode\ba_\lambda$ and $N\in\mode\ba_\mu$ be finite--dimensional and assume that \begin{equation} A/\tilde K^\lambda_M\tilde K^\lambda_N\cong
A/\tilde K^\lambda_M\oplus A/\tilde K^\lambda_N.\end{equation}
Then
\begin{equation}\label{gentensweyl}\bw^{\lambda+\mu}_A (M\otimes N)\cong \bw^\lambda_AM\otimes\bw^\mu_AN.\end{equation}
\end{cor}

\section{The algebra $\ba_\lambda$}
We continue to assume that $A$ is a finitely generated commutative associative algebra over $\bc$. Denote by $\max A$ the set of maximal ideals of $A$ and let ${\rm}J(A)$ be the Jacobson radical of $A$. In this section we shall identify the max spectrum of $\ba_\lambda$ and if $J(A)=0$ we shall also identify the algebra $\ba_\lambda$. As a consequence we also obtain a classification of the irreducible finite dimensional modules in $\cal I_\lambda^A$. Special cases of this classification were proved earlier in \cite{C1}, \cite{CP1} for $A=\bc[t,t^{-1}]$, in \cite{L} and \cite{R} in the case when $A$ is the polynomial ring in $k$ variables.

\subsection{}
  For  $r\in\bz_+$
 the symmetric group $S_r$ acts naturally on $A^{\otimes r}$ and $\max(A)^{\times r}$ and we let
 $(A^{\otimes r})^{S_r}$ be the corresponding ring of invariants and $\max(A)^{\times r}/S_r$ the set of orbits.
 If $r=r_1+\cdots +r_n$, then
  we regard $S_{r_1}\times \cdots\times S_{r_n}$ as a subgroup of $S_r$ in the canonical way, i.e $S_{r_1}$ permutes the first $r_1$ letters, $S_{r_2}$ the
  next $r_2$ letters and so on. Given $\lambda=\sum_{i\in I}r_i\omega_i\in P^+$, set
  \begin{gather}\label{defbba} r_\lambda=\sum_{i\in I} r_i,\ \ S_\lambda=S_{r_1}\times\cdots\times S_{r_n},\ \
  \ \ \mathbb{A}_\lambda=(A^{\otimes r_\lambda})^{S_\lambda},\\ \max(\mathbb A_\lambda)=(\max(A)^{r_1}/S_{r_1})\times\cdots\times(\max (A^{r_n})/S_{r_n}).\end{gather}
  The algebra $\mathbb A_\lambda$ is clearly finitely generated. For $\mathbb M\in \max \mathbb( A_\lambda)$, let $\ev_{\mathbb M}: \mathbb A_\lambda\to \bc$ be the corresponding algebra homomorphism.

We shall prove the following in the rest of the section.
  \begin{thm} \label{alambdaiso} \begin{enumerit}
  \item There exists a  canonical bijection $$\max\mathbb A_\lambda\to \max\ba_\lambda$$
  \item  Assume that  $\rm{J}(A)=0$ and let $\lambda\in P^+$. There exists an isomorphism  of algebras $$\tau_\lambda:\ba_\lambda\to\mathbb A_\lambda.$$
\end{enumerit}
  \end{thm}

\subsection{}
Let $\Xi$ be the monoid of finitely supported functions $\xi: \max(A)\to P^+$, where for $\xi,\xi'\in\Xi$ and $S\in\max A$, we define  \begin{gather*}
(\xi+\xi')(S)=\xi(S)+\xi'(S),\qquad \supp\xi=\{S\in\max(A): \xi(S)\ne 0\},
\qquad  \wt(\xi)=\sum_{S\in\max(A)}\xi(S).\end{gather*} Clearly $\wt:\Xi\to P^+$ is a morphism of monoids and we set
$$\Xi_\lambda=\{\xi\in\Xi:\wt\xi=\lambda\}.$$
Given $\xi\in\Xi_\lambda$, let $$K_\xi=\prod_{S\in\supp\xi}\!\!S,\ \quad \lie g_\xi= \lie g\otimes A/K_\xi,\ \quad  \bv_\xi= \bigotimes_{S\in\supp\xi} V(\xi(S)).$$ Then $\lie g_\xi$ is a finite--dimensional semi--simple Lie algebra and $\bv_\xi$ is an  irreducible finite--dimensional representation of $\lie g_\xi$ and hence of $\lie g\otimes A$ with action given by\begin{equation}\label{vxi}(x\otimes a)(v_1\otimes\cdots\otimes v_r)=\sum\limits_{k=1}^{r}\ev_{S_k}(a)(v_1\otimes \cdots \otimes xv_k\otimes\cdots\otimes v_r),\end{equation} where $S_1,\cdots, S_r$ is an enumeration of $\supp\xi$. Set $M_\xi= \br_A^\lambda\bv_\xi.$ By Lemma 5.2 we see that $\bv_\xi$ is the unique irreducible quotient of $\bw^\lambda_AM_\xi$ and hence  $$\bv_\xi\cong\bv^\lambda_AM_\xi.$$

Let $\lambda\in P^+$ and $M\in\rm{irr}\mode\ba_\lambda$. Since  $A/\tilde K^\lambda_M$ is  a finite--dimensional semi--simple algebra we know that $$\tilde K^\lambda_M= S_1\cdots S_r,\ \ \ \ r\in\bz_+,$$ where $S_1,\cdots, S_r$ are (uniquely defined up to permutation) maximal ideals in $A$. Moreover $\bv^\lambda_AM$ is a representation for the semi-simple Lie algebra $\lie g_M= \oplus_{k=1}^r\lie g\otimes A/S_i$. So there exist unique elements $\mu_1,\cdots,\mu_r\in P^+$ such that $$\bv^\lambda_AM\cong_{\lie g_M} V(\mu_1)\otimes \cdots\otimes V(\mu_r).$$   Define $\xi_M\in\Xi_\lambda$ by $$\xi_M(S_k)=\mu_k,\ \  1\le k\le r,\ \ \xi(S)=0,\ \ {\rm{otherwise}}. $$ Then $\bv^\lambda_AM\cong\bv_\xi$ as $\lie g\otimes A$--modules.
Summarizing, we have proved that:
\begin{prop} The assignment $\xi\to M_\xi$, (resp. $\xi\to \bv_\xi$) defines a natural bijection between $\Xi_\lambda$ and the set of isomorphism classes of irreducible representations of $\ba_\lambda$ (resp. isomorphism classes of  irreducible objects in $\cal I_A^\lambda)$. Moreover this bijection is compatible with the functor $\bv_A^\lambda$, in the sense that $$\bv_\xi\cong\bv^\lambda_AM_\xi.$$\hfill\qedsymbol
\end{prop}
Given $\xi\in\Xi_\lambda$, define $\ev_\xi:\bu(\lie h\otimes A)\to \bc$ by extending  $$\ev_\xi(h\otimes a)=\sum_{S\in\max A}\ev_S(a)\xi(S)(h).$$
\begin{cor}\label{anncont}  Let $\lambda\in P^+$. Then $$\ann_{\lie h\otimes A}w_\lambda\subset\bigcap_{\xi\in\Xi_\lambda}\ker\ev_\xi.$$
\end{cor}
\begin{pf} Let $u\in\bu(\lie h\otimes A)$ and assume that $uw_\lambda=0$. Since $\bv_\xi$ is a quotient of $W_A(\lambda)$ it follows that $u(\bv_\xi)_\lambda=0$. On the other hand it is clear from the definition of $\bv_\xi$ that $$(h\otimes a)(\bv_\xi)_\lambda=\ev_\xi(h\otimes a) (\bv_\xi)_\lambda,$$ and the corollary follows.
\end{pf}

\subsection{}
 The set $\Xi_\lambda$ also parametrizes the set $\max\mathbb A_\lambda$ as follows. Let  $\mathbb M\in\max(\mathbb A_\lambda)$ be the orbit of an element $(S_1,\cdots, S_{r_\lambda})\in\max(A)^{\times
r_\lambda}$. Define $\xi(\mathbb M)\in\Xi_\lambda$ by \begin{gather*}\xi(\mathbb M)(S)=\sum_{i\in I} p_i(S)\omega_i,\ \  S\in\max(A),\\
p_i(S)=\#\{p: \sum_{k=1}^{i-1}r_{k}< p\le \sum_{k=1}^{i}r_{k},\ \ \  S_p=S\}.\end{gather*} It is easily seen that the assignment $\mathbb M\to
\xi(\mathbb M)$ is well--defined bijection of sets $\max(\mathbb A_\lambda)\to \Xi_\lambda$ and part (i) of the Theorem is established.

\subsection{}
The algebra $\mathbb A_\lambda$ is generated by elements of the form  \begin{equation}\label{defnsym} \sym^i_\lambda(a)=1^{\otimes(r_1+\cdots r_{i-1})}\otimes \left(\sum_{k=0}^{r_i-1} 1^{\otimes k} \otimes a\otimes 1^{\otimes(r_i-k-1)}\right)\otimes 1^{\otimes (r_{i+1}+\cdots r_n)}, \ \ a\in A, \ i\in I.\end{equation} It is clear that the assignment $$\tilde{\tau}_\lambda(h_i\otimes a)=\sym^i_\lambda(a),\ \ i\in I, a\in A$$ extends to a surjective algebra homomorphism $\tilde{\tau}_\lambda:\bu(\lie h\otimes A)\mapsto\mathbb A_\lambda$. Moreover it is easily checked that \begin{equation}\label{kersa}\ev_{\xi(\mathbb M)}(h\otimes a)=\ev_{\mathbb M}\tilde\tau_\lambda(h\otimes a),\ \ h\in\lie h, \ a\in A.\end{equation}
\begin{lem} We have \begin{equation}\label{kers} \ker\tilde\tau_\lambda=\bigcap_{\mathbb M\in\max{\mathbb A_\lambda}}\ker\ev_{\mathbb M}\tilde{\tau_\lambda}=\bigcap_{\xi\in\Xi_\lambda}\ker\ev_\xi,\end{equation} and hence $\tilde\tau_\lambda$ induces a surjective homomorphism of algebras $\tau_\lambda:\ba_\lambda\to\mathbb A_\lambda$.
\end{lem}
\begin{pf} The first equality in \eqref{kers} is trivial since  $\rm{J}(\mathbb A_\lambda)=0$ if $\rm{J}(A)=0$. The second equality is immediate from \eqref{kersa} and the fact that $
\mathbb M\to \xi(\mathbb M)$ is bijective. The final statement of the Lemma is immediate from Corollary \ref{anncont}.
\end{pf}

\subsection{}
 It remains to prove that $\tau_\lambda$ is  injective. To do this we adapt an argument in \cite{FL}. Thus, we  identify a natural spanning set of $\ba_\lambda$ and prove that its image in  $\mathbb A_\lambda$ is a basis.  Fix an ordered  countable basis $\{a_r:r\in\bz_+\}$ of $A$ with $a_0=1$ and $a_r\in A_+$ for $r\ge 1$.

\begin{lem} The elements $$\{\prod_{i=1}^n\prod_{s=1}^{q_i}(h_i\otimes a_{i,s})w_\lambda: a_0< a_{i,1}\le\cdots\le a_{i,q_i},\ \  i\in I,\ \ q_i\le\lambda(h_i)\}$$ span $W_A(\lambda)_\lambda$.
\end{lem}
\begin{pf} It is clearly enough to prove that for each $i\in I$ and elements $1\le p_1\le\cdots\le p_\ell$,   $$\prod_{s=1}^\ell (h_i\otimes a_{p_s})w_\lambda\in{\rm{span}}\left\{\prod_{s=1}^{m}(h_i\otimes a_{r_s})w_\lambda: 1\le  r_1\le r_2\le\cdots \le r_{m},\ \ m\le \lambda(h_i)\right\}.$$ Since $$0=\prod_{s=1}^\ell(x_i^+\otimes a_{p_s})(x_i^-\otimes 1)^{\ell}= \prod_{s=1}^\ell (h_i\otimes a_{p_s})w_\lambda+ Hw_\lambda,\ \ \ell\ge \lambda(h_i)+1,$$ where $H$ is in the span of elements of the form $\prod_{s=1}^r (h_i\otimes a_{p_{j_s}})$ with $r<\ell$, the Lemma follows by a simple induction on $\ell$.
\end{pf}

\subsection{}
As a result of the Lemma we see that $\ba_\lambda$ is spanned by the image of the set $$\{\prod_{i=1}^n\prod_{s=1}^{m_i}(h_i\otimes a_{i,s}): a_{i,s}\in A_+,  a_{i,1}\le\cdots\le a_{i,m_i}, i\in I,\ \ m_i\le \lambda(h_i)\}.$$ The proof that  $\tau_\lambda$ is injective follows if we  prove that the set
$$\left\{\bigotimes_{s=1}^{m_1}\sym^1_\lambda(a_{1,s})\bigotimes\cdots \bigotimes_{s=1}^{m_n}\sym^n_\lambda(a_{n,s}):a_{i,s}\in A_+,  a_{i,1}\le\cdots\le a_{i,m_i}, i\in I,\ \ m_i\le \lambda(h_i)\right\}$$ is linearly independent in $\mathbb A_\lambda$. Since the tensor product of linearly independent sets is linearly independent it is enough to prove the following.
 Let $N\in\bz_+$ and for $b_1,\cdots, b_N\in A$ let $$\sym_N(b_1\otimes\cdots\otimes b_N)=\sum_{\sigma\in S_N}(b_{\sigma(1)}\otimes\cdots\otimes b_{\sigma(N)}).$$
 \begin{lem} The elements \begin{equation}\label{indep}\sym_N(a_{r_1}\otimes 1^{\otimes N-1})\sym_N(a_{r_2}\otimes 1^{\otimes N-1})\cdots\sym_N(a_{r_m}\otimes 1^{\otimes N-1}), 1\le r_1\le\cdots\le r_m,\ \ m\le N\end{equation} are linearly independent in $A^{\otimes N}$.
\end{lem}
\begin{pf} Set $$\mathbb U=\bigoplus_{0\le m\le N} A_+^{\otimes m}\otimes 1^{\otimes (N-m)},$$ and let $\bop: A^{\otimes N}\to \mathbb U$ be the canonical projection. The projection onto $\mathbb U$ of the elements in \eqref{indep} are
$$\sym_{r_m}(a_{r_1}\otimes a_{r_2}\otimes\cdots a_{r_m})\otimes 1^{N-m},\ \ 1\le r_1\le\cdots\le r_m,\ \ m\le N$$ and these are clearly linearly independent in $\mathbb U$ and the Lemma is proved.
\end{pf}

\section{The fundamental Weyl modules}
We use the notation of the previous sections freely.
Throughout this section we shall assume that $A$ is finitely generated.  Theorem \ref{alambdaiso}(i) applies and we have   bijections $\max(\ba_{\lambda}) \to \Xi_\lambda \to \max(\mathbb A_\lambda)$. Recall that $\max\mathbb A_\lambda$ is the set of orbits of the group $S_\lambda$ acting on $(\max A)^{\otimes r_\lambda}$. The orbits of maximal size (i.e those coming from an element of $(\max A)^{\otimes r_\lambda}$ with trivial stabilizer under the $S_\lambda$ action) correspond
under this bijection to the subset
$$\Xi_\lambda^{\rm{ns}}=\{\xi\in\Xi_\lambda: \xi(S)=\sum_{i\in I}m_i\omega_i,\ \ m_i\le 1\ \forall\ \ S\in\max A\}$$ of $\Xi$. The group $S_{r_\lambda}$ also acts on $(\max A)^{\otimes r_{\lambda}}$ by permutations and the orbits of this action can be naturally identified with a subset of $\max\mathbb A_\lambda$. The orbit of  points with trivial stabilizer under the $S_{r_{\lambda}}$ action  corresponds further
to the subset
$${}_1\Xi_\lambda^{\rm{ns}}=\{\xi\in\Xi_\lambda:  \xi(S)\in\{0,\omega_1,\cdots,\omega_n\},\ \ \forall\ \ S\in\max A\},$$ of $\Xi_\lambda^{ns}$.
Clearly $$\Xi_{\omega_i}={}_1\Xi^{ns}_{\omega_i}.$$

In this section we
shall analyze the modules $\bw^\lambda_AM_\xi$,  $\xi\in {}_1\Xi_\lambda^{\rm{ns}}$ when $\lie g$ is an algebra of classical type.  By Theorem~\ref{tensprod} we see that
\begin{equation}\label{tp1} \bw^\lambda_AM_\xi\cong\bigotimes_{S\in\supp\xi}\bw_A^{\xi(S)}M_{\xi_S},\ \ \supp\xi_S=\{S\},\ \ \xi_S(S)=\xi(S).
\end{equation}
This means that if $\xi\in {}_1\Xi^{\rm{ns}}_\lambda$, it is enough to analyze the modules $\bw^{\omega_i}_AM_\xi$, $i\in I$, $\xi\in\Xi_{\omega_i}$.

\subsection{}
 Assume from now on that $\lie g$ is of type $A_n$, $B_n$, $C_n$ or $D_n$. Assume also that the nodes of the Dynkin diagram of $\lie g$ are numbered as in \cite{B}.
Define a subset $J_0$ of $ I$  as follows: $$J_0=\begin{cases} I,\ \ \lie g \ \ {\rm{of\  type}}\  A_n,\  C_n,\\
\{n\},\ \ \lie g \ \ {\rm{of \ type}} \ B_n,\\
\{n-1,n\},\ \ \lie g \ \ {\rm{of \ type }}\ \ D_n.\end{cases}$$
Given $m,k\in\bz_+$, let $\boc(m)$ be the dimension of the space of  polynomials of degree $m$ in $k$--variables, i.e
$$\boc(m)=\#\{\bos=(s_1,\cdots, s_k)\in\bz_+^k: s_1+\cdots +s_k=m\}.$$
For $S\in\max A$ and $i\in I$ let  $\xi^i_S\in\Xi_{\omega_i}$ be given by requiring  $\supp\xi=S$.
\begin{thm} \label{fundamental-structure}  Assume that  $S\in\max A$ and that $\dim S/S^2=k$.
We have an isomorphism of $\lie g$--modules, \begin{gather}\label{ij0} \bw_{A}^{\omega_i}M_{\xi^i_S}\cong_{\lie g} V(\omega_i),\ \ i\in J_0,\\
\label{inotinj0}\bw_A^{\omega_i}M_{\xi^i_S}\cong_{\lie g}\bigoplus_{\{j: i-2j\ge 0\}} V(\omega_{i-2j})^{\oplus \boc(j)},\ \ i\notin J_0.\end{gather}
\end{thm}
\begin{rem} The theorem was proved  when $A$ is the polynomial ring in one variable in \cite{Cferm},\cite{CM}.  \end{rem}

\subsection{}
Before proving the theorem, we note  the following. Let $\dim_\lambda:\Xi_\lambda\to\bz_+$ be the function $\xi\to\dim\bw_A^\lambda M_\xi$.
\begin{cor} Let $A$ be a smooth irreducible algebraic variety.  The restriction of $\dim_\lambda$ to ${}_1\Xi^{{\rm{ns}}}_\lambda$ is constant.
\end{cor}
\begin{pf} Since $A$ is smooth and irreducible, it follows that $\dim S/S^2$ is independent of $S$ and hence by Theorem~\ref{fundamental-structure} we see that the corollary is true for $\omega_i$. The general case now follows from  \eqref{tp1}.
\end{pf}
\begin{rem}  In the special case when $A=\bc[t]$ the function $\dim_\lambda$ is constant on $\Xi_\lambda$. This was conjectured in \cite{CPweyl} and proved there for $\lie{sl}_2$. It was later proved in \cite{CL} for $\lie{sl}_{r+1}$, in \cite{FoL} for algebras of type $A,D,E$. The general case can be deduced by passing to the quantum group situation and using results in \cite{K}, \cite{BN}. No self--contained algebraic proof of this fact has been given for the non--simply laced algebras.

However, it is not true that if $A$ is an arbitrary  smooth irreducible variety, then  $\dim_\lambda$ is constant on $\Xi_\lambda^{\rm{ns}}$. As an example take $\lie g = \frak{sl}_3$, $A=\bc[t_1,t_2]$ and consider $\lambda=\omega_1+\omega_2$. Let $S,S'$ be the maximal ideals in $A$ corresponding to distinct points $(z_1,z_2,)$ and $(z_1',z_2')$.  Let $\xi, \xi'\in\Xi_\lambda$ be given by $$\xi(S)= \omega_1,\ \ \xi(S')= \omega_2,\ \ \xi'(S)=\lambda.$$
Then by Theorem \ref{fundamental-structure} $$\bw^\lambda_AM_\xi\cong_{\lie g} V(\omega_1)\otimes V(\omega_2),$$ and hence is nine--dimensional.

On the other hand the following argument proves that $\bw^\lambda_A(M_\xi')$ is at least 10--dimensional. Recall that  $V(\omega_1+\omega_2)\cong_{\lie g} \lie g_{\ad}$ where $\lie g_{\ad}$ is the adjoint representation of $\lie g$ and hence has dimension eight. Let $<\ ,\ >$ be the Killing form of $\lie g$.  A relatively straightforward check shows that if we set $W=\lie g_{\ad}\oplus\bc\oplus\bc$ and define an action of $\lie g\otimes A$ on $W$ by $$(x\otimes f)(y,z,z')= (f(z'_1,z'_2)[x,y],\  \frac{df}{dt_1}(z'_1,z'_2)<x,y>, \ \ \frac{df}{dt_2}(z'_1,z'_2)<x,y>),$$ then $W$ is a quotient of $\bw^\lambda_A(M_\xi')$.
\end{rem}

\subsection{}
The rest of the section is devoted to proving the theorem. We shall repeatedly use the following
\begin{equation}\label{haction}(\lie h\otimes S)(w_{\omega_i}\otimes M_{\xi^i_S})=0.\end{equation}
Given $\alpha\in R^+$, let $\varepsilon_i(\alpha)\in\{0,1,2\} $ be the coefficient of $\alpha_i$ in
 $\alpha$ and set $$\Ht\alpha= \sum_{j\in I}\varepsilon_j(\alpha),\ \qquad \lie n^-_r=\bigoplus_{\{\alpha\in R^+: \varepsilon_i(\alpha)=r\}}\lie g_{-\alpha}.$$ It is a simple matter to check that \begin{equation}\label{commutator}[\lie n^-_0,\lie n^-_0]=\lie n^-_0,\qquad [\lie n^-_0,\lie n^-_1]=\lie n^-_1,\qquad [\lie n_1^-,\lie n^-_1]=\lie n_2^-.\end{equation}
\begin{lem}\label{setofweights} We have
\begin{gather*}\left((\lie n^-_0\otimes A)\oplus(\lie n^-_1\otimes S)\oplus(\lie n^-_2\otimes S^2)\right)(w_{\omega_i}\otimes M_{\xi^i_S})=0.\end{gather*}
In particular $(\lie g\otimes S^2)\bw^{\omega_i}_AM_{\xi^i_S}=0$, i.e. $S^2\subset K^{\omega_i}_{M_{\xi^i_S}}$.
\end{lem}
\begin{pf}
It is trivial that  $$\lie n^+(x^-_j\otimes A)(w_{\omega_i}\otimes M_{\xi^i_S})=0,\ \ j\ne i,\qquad  \
   \lie n^+(x^-_i\otimes S)(w_{\omega_i}\otimes M_{\xi^i_S})=0 $$ Since $\omega_i-\alpha_j\notin P^+$ for all $i\in I$,
   it follows by elementary representation theory that
   $$(x^-_j\otimes A)(w_{\omega_i}\otimes M_{\xi^i_S})=0,\ \ j\ne i,\qquad  \  (x^-_i\otimes S)(w_{\omega_i}\otimes M_{\xi^i_S})=0.$$
Using \eqref{commutator} we see that   a straightforward induction on $\Ht\alpha$ proves the lemma.
\end{pf}

\subsection{}
We now prove by using  Lemma \ref{limitideal} and Lemma \ref{setofweights} that it suffices to prove Theorem \ref{fundamental-structure} in the case when $A$ is the polynomial ring in finitely many variables. For this, suppose that
 $B$ is a finitely generated algebra and let $S_B$ a maximal ideal in $B$. Let  $t_1, \ldots, t_k \in S$ be such that the images of these elements form a basis of $S_B/S_B^2$. Let $A = \bc[x_1, \ldots, x_k]$, and define an algebra homomorphism
$ A \longrightarrow B $ by extending the assignment $ x_i \mapsto t_i.$ Let $S_A$ be the ideal in $A$ generated by $x_1,\cdots,x_k$. Clearly $S_A$ maps to $S_B$ and we have a homomorphism of algebras $\phi: A/S_A^2\to B/S_B^2$.
Moreover, since $t_1,\cdots, t_k$ are linearly independent in $S_B/S_B^2$ it follows that $\phi$ is injective. Further, since $$\dim A/S_A^2=\dim B/S_B^2=k+1,$$ it follows that $\phi$ is an isomorphism of algebras.
We now have $$\bw_{B}^{\omega_i}M_{\xi^i_{S_B}} \cong \bw_{B/S_B^2}^{\omega_i}M_{\xi^i_{S_B}} \cong \bw_{A/(S_A)^2}^{\omega_i}M_{\xi^i_{S_A}}\cong \bw_A^{\omega_i}M_{\xi^i_{S_A}},$$ where the first and last isomorphisms follow from Lemma \ref{limitideal} and the isomorphism in the middle is induced by $\phi$.

\subsection{}
From now on we shall assume that $A=\bc[t_1,\dots, t_k]$ is the polynomial ring in $k$ variables. Moreover since the theorem is proved for $k =1$ in \cite{Cferm},\cite{CM}, we shall assume that $k>1$. In addition we may assume that $S$ is the maximal ideal generated by $t_1,\cdots, t_k$. There is no loss of generality in doing this for the following reason. Suppose that $S'$ is another maximal ideal corresponding to the point $\boz=(z_1,\cdots, ,z_k)\in\bc^k$. Consider the automorphism of $\phi_\boz: \lie g\otimes A\to\lie g\otimes A$ given by $x\otimes t_r\to x\otimes (t_r-z_r)$, $x\in\lie g$, $1\le r\le k$. It is not hard to check that $$\bw^{\omega_i}_AM_{\xi^i_S}\cong\phi_\boz^*\bw^{\omega_i}_AM_{\xi^i_{S'}}.$$

\subsection{}
Let $A_+$ be the subspace of  polynomials with constant term zero. Since $\lie g\otimes A_+$ is an ideal in $\lie g\otimes A$, to prove \eqref{ij0} it suffices to show that for all $\alpha\in R^+$ and $a\in A_+$,\begin{equation}\label{gmap} (x^-_\alpha\otimes a)(w_{\omega_i}\otimes M_{\xi^i_S})\in\bu(\lie g)(w_{\omega_i}\otimes M_{\xi^i_S}).\end{equation}
Let $C = \bc[t]$, where $t$ is an indeterminante. Consider the map $\lie g\otimes C\to\lie g\otimes A$ given by $x\otimes t\to x\otimes a$. By Proposition {\ref{canonical} there exists a map of $\lie g\otimes C$--modules $\bw^{\omega_i}_CM_{\xi^i_S}\to \bw^{\omega_i}_AM_{\xi^i_S}$. Since the theorem is known for $C$, it follows that $$(x_\alpha^-\otimes t)(w_{\omega_i}\otimes M_{\xi^i_S})\in\bu(\lie g)(w_{\omega_i}\otimes M_{\xi^i_S})\subset\bw^{\omega_i}_CM_{\xi^i_S},$$ which proves \eqref{gmap}.

\subsection{}
The rest of the section is devoted to proving \eqref{inotinj0} and hence we may and will  assume  that $\lie g$ is of type $B_n$ or $D_n$. For $j\in I$, $j\ge 2$, set  $ \omega_j - \omega_{j-2}=\theta_j$. Then one checks easily \cite{Humphreys}  $$\theta_j\in
R^+,\ \quad \theta_{j-2}-\theta_j=\alpha_{j-3}+2\alpha_{j-2}+\alpha_{j-1}, \quad \theta_j-\alpha_r\in R^+\ \ \iff r=j.$$ where we understand that $\alpha_{-1}=0$.
\begin{prop}\label{quot1} Let $i\in I$, $1\le \ell,m\le k$, and set $v_\ell=(x^-_{\theta_i}\otimes t_\ell )(w_{\omega_i}\otimes M_{\xi_S^{i}})$.
 Then\begin{gather}\label{quot}(\lie n^+\otimes A) v_\ell=0,\ \ \qquad  (\lie h\otimes S)v_\ell=0.\end{gather} In particular, the $\lie g\otimes A$--submodule of $\bw^{\omega_i}_AM_{\xi^i_S}$ generated by $v_\ell$ is  a quotient of $\bw^{\omega_{i-2}}_AM_{\xi^{i-2}_S}$.
 Further, we have \begin{gather} \label{reorder}(x^-_{\theta_{i-2}}\otimes t_m)v_\ell=(x^-_{\theta_{i-2}}\otimes t_\ell)v_m.\end{gather}
\end{prop}
\begin{pf} Note that  $(\lie n^+\otimes
S)v_\ell$ and $ (\lie h\otimes S)v_\ell$ are both contained in $(\lie g\otimes S^2)(w_{\omega_i}\otimes M_{\xi^i_S}) $ and hence by Lemma \ref{setofweights}$$(\lie n^+\otimes
S)v_\ell\ \ =0\ \ =\  (\lie h\otimes S)v_\ell.$$
Since $S$ is maximal, \eqref{quot} follows if we  prove that $$(\lie n^+\otimes 1)
v_\ell=0.$$ Since $$[x^+_{j}, x^-_{\theta_i}]=0,\ j\ne i, \ \ {\rm{and}}\
\  \varepsilon_i(\theta_i-\alpha_i)=1,$$ we see that Lemma
\ref{setofweights} gives $(x_j^+\otimes 1)v_\ell=0$ for all $j\in
I$. The second statement of the proposition is now clear.
Hence we have  by Lemma \ref{setofweights} that
 $$(x^-_\alpha\otimes S)v_\ell=0\ \qquad
 {\rm{if}}\ \ \varepsilon_{i-2}(\alpha)\ne 2. $$
Writing $$x^-_{\theta_{i-2}}=[x^-_{i-2},
[x^-_{\alpha_{i-3}+\alpha_{i-2}+\alpha_{i-1}}, x^-_{\theta_i}]],$$
and using Lemma \ref{setofweights} we get
\begin{gather*}(x^-_{\theta_{i-2}}\otimes
t_m)v_\ell=(x^-_{\theta_{i-2}}\otimes t_m)(x^-_{\theta_{i}}\otimes
t_\ell)(w_{\omega_i}\otimes M_{\xi^i_S}) =(x^-_{\theta_{i}}\otimes
t_\ell)(x^-_{\theta_{i-2}}\otimes t_m)(w_{\omega_i}\otimes M_{\xi^i_S})\\=
(x^-_{\theta_{i}}\otimes
t_\ell)x^-_{i-2}x^-_{\alpha_{i-3}+\alpha_{i-2}+\alpha_{i-1}}(
x^-_{\theta_i}\otimes t_m)(w_{\omega_i}\otimes M_{\xi^i_S})\\ =
(x^-_{\theta_{i-2}}\otimes t_\ell)v_m+ X(w_{\omega_i}\otimes M_{\xi^i_S}),\end{gather*} where $X$
is a linear combination of the elements
\begin{gather*}x^-_{i-2}x^-_{\alpha_{i-3}+\alpha_{i-2}+\alpha_{i-1}}(x^-_{\theta_{i}}\otimes
t_\ell)(x^-_{\theta_i}\otimes t_m),\ \
x_{i-2}^-(x^-_{\theta_i+\alpha_{i-3}+\alpha_{i-2}+\alpha_{i-1}}\otimes
t_\ell)(x^-_{\theta_i}\otimes t_m),\\
x^-_{\alpha_{i-3}+\alpha_{i-2}+\alpha_{i-1}}(x^-_{\theta_i+\alpha_{i-2}}\otimes
t_\ell)(x^-_{\theta_i}\otimes t_m).\end{gather*}
But by Lemma~\ref{setofweights} all these terms act as zero on $(w_{\omega_i}\otimes M_{\xi^i_S})$, since
$(x^-_{\theta_i}\otimes t_m)(w_{\omega_i}\otimes M_{\xi^i_S})$ generates a quotient of
$\bw_{A}^{\omega_{i-2}}M_{\xi_S^{i-2}}$ and
$$\varepsilon_{i-2}(\theta_i+\alpha_{i-2}+\alpha_{i-1}+\alpha_{i-3})=1=\varepsilon_{i-2}(\theta_i+\alpha_{i-2}).$$
\end{pf}
The following is now immediate.
\begin{cor} Given $i,\ell\in I$ with $2\ell\le i$,   and  $  r_s\in\{1,\cdots, k\}$, $1\le s\le \ell$, the elements,
$$v(r_1,\cdots ,r_\ell)=(x^-_{\theta_{i -2\ell}} \otimes
t_{r_{\ell}})\cdots (x^-_{\theta_{i - 2}} \otimes
t_{r_{2}})(x^-_{\theta_{i}} \otimes t_{r_1}).(w_{\omega_i}\otimes M_{\xi^i_S})$$
generate a submodule of $\bw_{A}^{\omega_i}M_{\xi^{i}_S}$ which is a quotient
of $\bw_{A}^{\omega_{i-2\ell}}M_{\xi_S^{i-2\ell}}$. Moreover if $\sigma\in S_\ell$,
we have,$$v(r_1,\cdots ,r_\ell) =v(r_{\sigma(1)},\cdots
,r_{\sigma(\ell}).$$
\end{cor}

\subsection{}
Suppose that $\alpha\in R^+$ is such that
$\varepsilon_i(\alpha)=2$. Then we can write
$\alpha=\gamma + \beta+\theta_i$ for some $\beta,\gamma\in R^+$ with
$\varepsilon_i(\beta)=\varepsilon_i(\gamma)=0$. This implies that
$x^-_\alpha=c[x^-_\beta,[x^-_\gamma,x^-_{\theta_i}]],$ for some
non--zero $c\in\bc$ and hence  $$(x^-_\alpha\otimes t_\ell)(w_{\omega_i}\otimes M_{\xi^i_S}) = c[x^-_\beta,[x^-_\gamma,x^-_{\theta_i}\otimes t_\ell]]
(w_{\omega_i}\otimes M_{\xi^i_S})\in\bu(\lie g)(x^-_{\theta_i}\otimes
t_\ell)(w_{\omega_i}\otimes  M_{\xi^i_S}).$$ Proposition \ref{quot1} now gives,
 $$\bw_{A}^{\omega_i}M_{\xi^i_S}= \bu(\lie g)(w_{\omega_i}\otimes M_{\xi^i_S})\oplus\sum_{\ell=1}^k\bu(\lie g\otimes A)(x^-_{\theta_i}\otimes t_\ell)(w_{\omega_i}\otimes M_{\xi^i_S}) $$ as $\lie g$--modules.
 Using Corollary \ref{quot1} we find
$$\bw_{A}^{\omega_i}M_{\xi^i_S}= \bu(\lie g)(w_{\omega_i}\otimes M_{\xi_S^i})\bigoplus_{0\le
2l\le i} \; \; \left(\sum_{0\le r_1\le\cdots\le r_\ell}\bu(\lie
g)v(r_1,\cdots,r_\ell)\right),$$ which proves that $$\Hom_{\lie g}(V(\mu), \bw^{\omega_i}_AM_{\xi_S^{i}})=0,\ \ \mu\ne i-2j,\ \  \dim\Hom_{\lie g}(V(\omega_{i-2j}), \bw^{\omega_i}_AM_{\xi_S^{i}})\le \boc(j).$$

\subsection{}
To complete the proof  it suffices to prove that the elements $v(r_1,\cdots,r_l)$ are linearly independent for all $i,\ell\in I$ with $2\ell\le i$ and  $  r_s\in\{1,\cdots, k\}$, $1\le s\le \ell$. We do this as in \cite{CM} by explicitly constructing a module which is a quotient of $\bw^{\omega_i}_AM_{\xi_S^{i}}$ and where these elements are linearly independent.
Suppose that $V_s$ for $0\le s\le\ell$ are $\lie g$--modules such that
\begin{equation}\label{cons} \Hom_{\lie g}(\lie g\otimes V_s, V_{s+1})\ne 0, \ \ \Hom_{\lie g}(\wedge^2(\lie g)\otimes V_s,V_{s+1})=0.\end{equation}  Set $V=\oplus_{s=0}^\ell V_s$ and fix non--zero elements $p_s\in \Hom_{\lie g}(\lie g\otimes V_s, V_{s+1})$ for $0\le s\le k$. Define a $\lie g\otimes A$--module structure on $V\otimes A$ by:
  \begin{gather*}(x\otimes 1)(v\otimes a)= xv\otimes a,\ \ (x\otimes t_r)(v\otimes a)=p_s(x\otimes v)\otimes at_r,\ \ x\in\lie g,\ \ a\in A \ \ 1\le r\le k,\\
   (x\otimes S^2)(v\otimes a)=0.\end{gather*}
   To see that this is  an action, the only non--trivial part is to notice that, \begin{gather*}[x\otimes t_r, y\otimes t_m](v\otimes c)=p_{s+1}(x\otimes p_s(y\otimes v))\otimes t_rt_m c-p_{s+1}(y\otimes p_s(x\otimes v))\otimes t_rt_m c,\\=
p_{s+1}(p_s\otimes 1)((x\otimes y-y\otimes x)\otimes v)\otimes t_rt_\ell c=0,\end{gather*} where the last equality follows by noticing that  $p_{s+1}(p_s\otimes 1)\in\Hom_{\lie g}(\lie g\otimes\lie g\otimes V_s,V_{s+1})$ and using \eqref{cons}. 

It was shown in \cite{CM} that the modules $V(\omega_{i-2s})$, $0\le 2s\le i$ satisfy \eqref{cons} and also that
  $$p_s(x^-_{\theta_{i-2s-2}}\otimes v_{\omega_{i-2s}})=v_{\omega_{i-2s-2}}.$$
 and hence we can apply the preceding construction to this family of modules. Consider the $\bu(\lie g\otimes A)$--module $\bar W$ generated by $v_{\omega_i}\otimes 1$. It is clear that $$(\lie n^+\otimes A)(v_{\omega_i}\otimes 1)\ = \ 0\ = (\lie h\otimes S)(v_{\omega_i}\otimes 1),$$ since $\omega_{i-2}<\omega_i.$  Hence $\bar W$ is a quotient of $\bw^{\omega_i}_AM_{\xi_S^{i}}$.  Moreover, it is simple to check now that
$$(x^-_{\theta_{i -2\ell}} \otimes t_{r_{\ell}})\cdots (x^-_{\theta_{i - 2}} \otimes t_{r_{2}})(x^-_{\theta_{i}} \otimes t_{r_1}).v_{\omega_i}=v_{\omega_{i-2\ell}}\otimes t_{r_1}\cdots t_{r_\ell}.$$ Since these elements are manifestly linearly independent the result follows.

\end{document}